\documentclass[12pt]{article}

\usepackage[titletoc,title]{appendix}

 \newcommand{\lab}[1]{\label{#1}}

\usepackage{graphics,graphicx}
\usepackage{amsmath, amssymb}

\def\fullpage
{
\addtolength{\topmargin}{-1.5 cm}
\addtolength{\oddsidemargin}{-1.5 cm}
\addtolength{\textwidth}{+3 cm}
\addtolength{\textheight}{+3 cm}
}
\fullpage

\catcode`@=11
\@addtoreset{equation}{section}

\catcode`@=12

\def\non{\nonumber}
\newcommand{\be}{\begin{equation}}
\newcommand{\ee}{\end{equation}}
\newcommand{\eq}{\begin{equation}}
\newcommand{\en}{\end{equation}}
\newcommand{\ignore}[1]{}
\def\blackslug{\hbox{\kern1pt\vrule height6pt width4pt  depth1pt\kern1pt}}
\def\qed{\penalty 500\hbox{\quad\blackslug}\ifmmode\else\par
         \vskip4.5pt plus3pt minus2pt\fi}
\newtheorem{thm}{Theorem}[section]
\newtheorem{lemma}[thm]{Lemma}
\newtheorem{cor}[thm]{Corollary}
\newtheorem{prop}[thm]{Proposition}

\def\proof{\par\noindent{\bf Proof.\enspace}\rm}

\def\be{\begin{equation}}
\def\ee{\end{equation}}
\def\bea{\begin{eqnarray}}
\def\eea{\end{eqnarray}}
\def\bean{\begin{eqnarray*}}
\def\eean{\end{eqnarray*}}
\newcommand{\eqn}[1]{(\ref{#1})}
\newcommand{\bel}[1]{\be\lab{#1}}

\def\BP{{\bf P}}

\def\eps{\epsilon}
\def\epshat{\overline{\eps}}
\def\error{n^{-\epshat}}

\def\d{\delta}
\def\ga{\gamma}
\def\la{\lambda}
\def\m{\mu}
\def\n{\nu}

\def\tilj{{\tilde j}}
\def\tilbfg{\tilde{\bf g}}
\def\bfg{{\bf g}}
\def\bfgt{{\bf g_t}}
\def\EKn{\Omega}
\def\smalls{{\cal S}}

\def\larges{{\cal L}}
\def\types{{\cal T}}
\def\type{\tau}
\def\hamma{\hat{\gamma}}
\def\hatma{\overline{\gamma}}
\def\clusters{{\cal K}}
\def\bz{{\bf 0}}
\def\Gnp{{\cal G}(n,p)}
\def\pr{{\bf P}}
\def\ex{{\bf E}}

\def\vmax{{v_{\rm max}}}
\def\vmaxx{{v_{\rm max}^{(1)}}}
\def\ut{{\frac{\lambda_u}{\lambda_t}}}
\def\functions{{\cal F}}
\def\smallfunctions{{\cal F}_\smalls}
\def\hfunctions{{\cal F}_1}
\def\unavoid{{\cal U}}
\def\class{{\cal C}}
\def\extensions{{\cal E}}
\def\Rset{{\cal R}}
\def\subs{\subseteq}

\def\stcolon{\; : \:}
\def\mean{\lambda}

\def\calP{{\cal P}}

\def\Gnp{{\cal G}(n,p)}
\def\Gnm{{\cal G}(n,m)}

\newcommand{\real}{\ensuremath {\mathbb R} }

\begin{document}

\title{The probability of nonexistence of a subgraph in a moderately sparse random graph }
\author{ Dudley Stark\thanks{Research  initially undertaken  while this author was a member of the Department
of Mathematics and Statistics,
University of Melbourne}\\
School of Mathematical Sciences\\
 Queen Mary College\\ University of London
 \and  Nick Wormald\thanks{Research supported by the Australian Laureate Fellowships grant FL120100125. Also supported by the
Australian Research Council  while this author was a member of the Department
of Mathematics and Statistics,
University of Melbourne, and by the Canada Research Chairs program and NSERC while in the Department of Combinatorics and Optimization, University of Waterloo} \\
School of Mathematical Sciences\\
 Monash University}

\date{}
\maketitle

\begin{abstract}
We develop  a general procedure that finds recursions for 
statistics counting isomorphic copies of a graph $G_0$  in the common random graph models $\Gnm$ and $\Gnp$. Our results apply 
 when the average degrees  of the random graphs  are 
below  the  threshold at which each edge is included in a copy of 
 $G_0$.  This extends an argument given earlier by the second author for $G_0=K_3$ with a more restricted range of average degree.
For all strictly balanced subgraphs $G_0$, our results give much information on the
  distribution of the number of copies of $G_0$ that are not in large ``clusters'' of copies. The probability that a random graph in
$\Gnp$ has no copies of $G_0$ is shown to be given asymptotically by  the exponential of a power
series in $n$ and $p$, over a fairly wide range of
$p$. A corresponding result is also given for $\Gnm$, which gives an asymptotic formula for the number of graphs with $n$ vertices, $m$ edges and no copies of $G_0$, for the applicable range of $m$. An example is given, computing the asymptotic probability that a random
graph has no triangles for
$p=o(n^{-7/11})$ in
${\cal G}(n,p)$ and for $m=o(n^{15/11})$ in ${\cal G}(n,m)$,
extending results of the second author. 
\end{abstract}

\newpage

\section{Introduction}

Our topic is the number of subgraphs of a random graph that are isomorphic to some given graph $G_0$. The perturbation method of~\cite{W} is used to derive
recursions of ratios of random graph statistics describing
the occurence of different types of clusters formed as edge-overlapping
groups of copies of $G_0$.
These recursions are used to
investigate the probability of no occurrences of $G_0$,
as well as other aspects of the
distribution of clusters.
For certain graphs $G_0$ and restrictions on $p$, we show that the
probability that there are no copies of the graph in $\Gnp$ is
the exponential an appropriate truncation of a power series
in $n$ and $p$, with error factor $(1+o(1))$. (As is usual, $\Gnp$ denotes the random graph on $n$ vertices
obtained by choosing
each edge in the graph to be present independently with probability $p$
and $\Gnm$ denotes the random graph on $n$ vertices obtained by
choosing uniformly at random from the $\binom{\binom{n}{2}}{m}$
graphs having $m$ edges.)
By considering recursions involving both  $G_0$ and
isolated edges, we build on this result to show that
the probability that there are no copies of $G_0$ in $\Gnm$ is
given in the same way but by a different power series in $n$ and $d$, where
\eq\lab{ddef}
d=\frac{m}{\binom{n}{2}},
\en
under corresponding restrictions on $d$.

Let $\nu(G)$ and $\mu(G)$ denote the number of vertices and number
of edges of a graph $G$. 
A graph
$G_0$ is {\em strictly balanced} if all its subgraphs are strictly less dense than $G_0$; that is,
$$
\frac{\m(G_0)}{\n(G_0)} > \frac{\m(G_1)}{\n(G_1)}
$$
for all nontrivial proper subgraphs $G_1$ of $G_0$. For example, the graph $K_n$ is
strictly balanced for all $n\ge 2$, as is every cycle.
Let $G_0$ be strictly balanced, and let $X$
be the number of copies of $G_0$
in the random graph $\Gnp$.
Let $\chi>0$ be defined by
\be\lab{chidef}
\chi=\chi(G_0)=\max_{G_1\in\extensions}\frac{\nu(G_0)-\nu(G_1)}{\mu(G_0)-
\mu(G_1)}.
\ee
We will restrict the growth of $p$ to 
$p=O(n^{-\chi-\epsilon})$ for some $\eps>0$.
The reason for this restriction
is that when $p$ is a little larger than $n^{-\chi}$   (sometimes called the {\em 2-threshold}), 
 each edge of   $\Gnp$ will expect to be contained in many copies of $G_0$. Thus, there will be subgraphs consisting of arbitrarily large numbers of copies of  $G_0$ ``chained'' together by shared edges. In this   case our analysis will not apply, since it relies on a copy of $G_0$ being unlikely to overlap with any others, as happens when
 restricting to $p=O(n^{-\chi-\epsilon})$.

Here is our main result. Note that $\chi$ should not be confused with the chromatic number, which does not appear in this paper.
\begin{thm}\lab{t:main}
Let $G_0$ be strictly balanced and put $\chi=\chi(G)$. Let $X$ be the  
number of copies
of $G_0$ in $\Gnp$, or let  $X$ be the number of copies
of $G_0$ in $\Gnm$ and set  $p= m/\binom{n}{2}$.
In each case, there is a formal power series $F= F(G_0) =\sum_{\ell\geq 0}c_{\ell}n^{i_\ell}p^{j_\ell}$, with   $i_\ell$ and $j_\ell$ strictly positive
for all
$\ell$, depending only on $G_0$, such that the following holds. For any $\eps>0$, if
  $p=O(n^{-\chi-\eps})$, then
\eq\lab{maineq}
\BP(X=0)=\exp\left(\sum_{\ell =0}^{M_\eps}c_{\ell }n^{i_\ell}p^{j_\ell} +o(1)\right),
\en
where the bound implicit in $o(1)$ is uniform over all such $p$ (but depends on $\epsilon$), and $M_\eps$ is a constant depending only on $\eps$
and $G_0$. Moreover, $\ell>M_\eps$ if and only if  $i_\ell<j_\ell(\chi+\eps)$.  
 \end{thm}

\noindent
{\bf Remarks}
\smallskip

\noindent
1. 
The theorem immediately  gives an asymptotic formula for the number of $G_0$-free graphs on $n$ vertices and $m$ edges, for the values of $m$ covered, by multiplying the $\Gnm$ case of~\eqn{maineq} by
${n(n-1)/2 \choose m}$.
\smallskip

\noindent
2.
Note that $i_\ell<j_\ell(\chi+\eps)$ if and only if  $n^{i_\ell}p^{j_\ell}= o(1)$ when $p=n^{-\chi-\epsilon}$, so each term with $\ell>M_\eps$ is $o(1)$. We also note that the issue of non-convergence of the  power series $F(G_0)$ for a given fixed $n$ and $p$   is not relevant in the present context.
  \smallskip

\noindent
3.
The proof of the theorem contains a definition of the  coefficients $c_\ell$ in Theorem~\ref{t:main} in terms of an algorithm by which they may be computed. It involves summing over a set of graphs whose size is  bounded for fixed $\eps>0$, but not as $\eps\to 0$. 
\smallskip

We next give two specific examples of the main result, by restricting to that case that $G_0$ is a triangle, or $K_3$, and   computing only the first few terms of the power series explicitly. 
\begin{thm}\lab{t:Gnp}
 If $p=p(n)=o(n^{-7/11})$,  the probability that the random graph $\Gnp$ is triangle-free is asymptotic to
$$
\exp\left(
-\frac{1}{6}n^3p^3 + \frac{1}{4}n^4p^5 - \frac{7}{12}n^5p^7
+\frac{1}{2}n^2p^3 -\frac{3}{8}n^4p^6 + \frac{27}{16}n^6p^9
\right).
$$
\end{thm}
Similarly, we determine the coefficients
$c_\ell$
  in the case of $\Gnm$ where
$G_0=K_3$ and $d=o(n^{-7/11})$,
 or equivalently  $m=o(n^{15/11})$,
in  the next theorem.
\begin{thm}\lab{t:Gnm}
 If $m=m(n)=o(n^{15/11})$, the probability that the random graph $\Gnm$ is triangle-free is asymptotic to 
$$
 \exp\left(
-\frac{1}{6}n^3d^3 
-\frac{1}{8}n^4d^6
\right),
$$
where $d=m/\binom{n}{2}$.
\end{thm}
These two results on triangles agree with and extend those of the second author in~\cite{W}, which applied for $p=o(n^{-2/3})$ and extended earlier results of Frieze~\cite{F}.

 For $\Gnm$, the expected value of $X$   is easily found to be 
\bean
  \la(G_0)&:= &
 \binom{n}{\nu } 
\binom{m}{\mu }
\binom{\binom{n}{2}}{\mu }^{-1}   \nu!
|\mbox{aut}(G_0)|^{-1}\\
&\sim& \hat \la(G_0):= \frac{(2m)^{ \mu }}{n^{2\mu-\nu} |\mbox{aut}(G_0)|}
\eean
where $\nu=\nu(G_0)$, $\mu=\mu(G_0)$, and  $|\mbox{aut}(G_0)|$ denotes the number of automorphisms of $G_0$. Ruci{\'n}ski~\cite{R} showed that the distibution of $X$ is asymptotically Poisson essentially for $d$ up to $n^{-\chi}$.
Frieze~\cite[Remark~2, P.69]{F} raised the possibility that, for the same range of $d$, the number of graphs with $k$ copies of $G_0$ in $\Gnm$ is asymptotic to the  probability that the Posson random variable with mean $\hat \la(G_0)$ is equal to $k$, for all ``small'' $k$. Theorem~\ref{t:Gnm} shows (for the first time!) that this is false   in particular for $k=0$ and $G_0=K_3$, since in this case, $\chi=1/2$ but already for $m=n^{4/3}$, other terms are entering the asymptotic formula in a significant way. Moreover, the situation is not remedied by using (the more natural) Poisson with mean $\la(G_0)$, since 
$\la(G_0) =  \frac{1}{6}n^3d^3 -\frac{1}{2}n^2d^3 +o(1)$ (using $nd^2=O(m^2/n^3)=o(1)$ for the range of $m$ under consideration).
  
We note that it may be possible to modify our approach to cater also for subgraphs that are not strictly balanced. In some cases, for instance where $G_0$ has a unique densest subgraph, the desired result can be deduced immediately from our results. However, other cases are more delicate, with different subgraphs of $G_0$ `competing'. One would need to incorporate considerations similar to those in the determination the threshold of appearance of $G_0$, as was done by Bollob{\'a}s~\cite{B}. 
 
Our concern 
here is to obtain an asymptotic formula for the probability that a random graph in $\Gnp$ or $\Gnm$ is $G_0$-free, for  a fixed graph $G_0$, where the density of the random graph is small enough that there are no large clusters of copies of $G_0$. 
Our methods will not work for the denser case, but some results are already known there, and for arbitrary densities. Recall, as in Remark~2 above, that for  $\Gnm$ our problem is equivalent to enumerating $m$-edged graphs with a forbidden subgraph. The classic paper of Erd{\H o}s, Kleitman and Rothschild~\cite{EKR}   gives the number of triangle-free graphs with $n$ vertices, in total, asymptotically (and asymptotics of the logarithm of the number when $G_0=K_t$). These results also demonstrate the connection between enumeration and the extremal numbers of edges for $G_0$-free graphs. There are many other similar results, which we refrain from mentioning as they do not take into accoung the edge density of the host graph.  More related to   the problem at hand, Pr{\"o}mel and Steger~\cite{PS96} found an asymptotic formula for the number of triangle-free graphs with 
$n$ vertices and $m$ edges when   $m> c n^{7/4} \log n$, by showing that they are almost all bipartite. This was extended by D. Osthus, H.J.~Pr\"omel and A. Taraz~\cite{OPT} to cover all $m$ that are at least slightly above $n^{3/2}$.  
Before this, {\L}uczak~\cite{L} had found asymptotics of the logarithm of the number.

For more general subgraphs than the triangle, and general $p$,  asymptotic formulae for the actual numbers (or probabilities) are elusive. The {\em logarithm} of the probability that $\Gnp$ is  $G_0$-free was estimated within a constant factor by Janson, \L uczak and Ruci\'nski~\cite{JLRinequal}. This was extended by Pr{\"o}mel and Steger~\cite{PS} to similar bounds on $\pr\big(\mbox{$\Gnm$ is  $G_0$-free}\big)$.
 
Many   results are known on the distribution of the number of copies of a fixed subgraph in $\Gnp$ and $\Gnm$; see for example~\cite[Chapter 6]{JLR}, but this is not our concern in this paper.

Our basic approach, and its background, are discussed in~\cite{W}.  
 The proof for $\Gnp$ estimates ratios of numbers of graphs using induction on the numbers of edge-overlapping clusters of copies of $G_0$ up to a given size; for $\Gnm$ the number of edges not in copies of $G_0$ is also used,  and the base step of this induction is essentially given by the $n$-vertex graph with no edges. 
There are two major extensions to the argument in~\cite{W}. One is that the graph $G_0$ is no longer restricted to $K_3$. This extension requires mainly graph theoretic arguments related to the ways that multiple copies of a graph can overlap.  The other is that the range of $p$ permits  edge-overlapping clusters containing arbitrarily many copies of $G_0$ to appear in the typical random graph under consideration. Thus our asymptotic estimates involve polynomials of unbounded size, and this poses significant problems in characterising and managing those estimates  (see Corollary~\ref{expanded} for example). 

The working assumption on $p=p(n)$ we will make in our proofs is
$p=n^{-\kappa+o(1)}$ where $\kappa\geq\chi+\epsilon$ is {\em fixed}. 
This assumption can be weakened to obtain asymptotic results that hold uniformly over more general $p=p(n)=O(n^{-\chi-\epsilon})$ 
by using
the following lemma. Here $a$ and $b$ are finite but the same result holds (with appropriate interpretation) without this assumption.    
\begin{lemma}\lab{panything}
For a closed interval $ [a,b] $, suppose that $f(n,p)$ is a function
such that $f(n,p)\to 0$ as $n\to\infty$
for all $p$ of the form $p=n^{-\kappa+o(1)}$ when $\kappa\in [a,b]$ is fixed.
Then
for fixed $\epsilon>0$,
$f(n,p)\to 0$ uniformly  for all $p(n)$ satisfying $p(n)= n^{-\kappa(n)}$
with $\kappa(n)\in[a,b]$ for all $n$.
\end{lemma}
\proof
 If  $p(n)$ satisfies $-\log_n p \in[a,b]$ for all $n$, then any subsequence of  $\big(p(n)\big)_{n\ge 1}$   has a subsubsequence for which $-\log_n p\to \kappa'$ for some fixed  $\kappa'\in[a,b]$. On this subsubsequence, $f(n,p)\to 0$ by assumption. So the lemma follows from the subsubsequence principle (see~\cite[p.12]{JLR}) applied to the sequence $\big(f(n,p(n))\big)_{n\ge 1}$ .  \qed

Our results will give information on the distribution of the number of copies of a strictly balanced subgraph, not just the probability that the number is 0, but we postpone this investigation to another paper.   We believe that it should be possible to modify our approach so as to obtain   accuracy in the formulae to any desired power of $n^{- {1}}$. Specifically, the power series in Theorem~\ref{t:main} should give valid lower order correction terms to the asymptotic formulae.    However, we have avoided attempting  this and there are some steps in the present argument that would have to be replaced in order to carry it out.
 
 Some basic definitions are made and results are proved in Section~\ref{s:recursions}; the $\Gnp$ case of Theorem~\ref{t:main} is proved in 
in Section~\ref{s:main};
the $\Gnm$ case is proved in
in Section~\ref{s:maingnm};
 Theorems~\ref{t:Gnp} and~\ref{t:Gnm} 
are proved in Appendix~\ref{s:tri}.

\section{
Clusters and recursions for counting maximal clusters}\lab{s:recursions}
 
We assume for a general framework that  $\EKn$ is any finite set. A family
$\clusters$ of subsets of  $\EKn$ is called a {\em clustering} if
      $C_1\in \clusters$, $C_2\in \clusters$ and $C_1\cap C_2\ne\emptyset$
imply that
$C_1\cup C_2\in\clusters$. The elements of $\clusters$ are called {\em
clusters}.

We will consider here only
the case that $\EKn = \EKn_n$ is the set of edges of the complete graph
$K_n$ on $n$ vertices, although the same
principles can also be applied to clusterings in general. As a further
restriction, to focus on small subgraph counts, we only consider very special
clusterings, for which  simplification occurs by taking
advantage of the symmetries of $K_n$. We take a fixed
graph $G_0$ throughout this paper, and will investigate the distribution of
the number of subgraphs of a random graph isomorphic to $G_0$. The edge set of any
subgraph of $K_n$ isomorphic to $G_0$ is called an {\em elementary $G_0$-cluster}. 
Mostly, we deal with
the minimal clustering which has every
elementary $G_0$-cluster as a member. We call this the {\em
$G_0$-clustering} of $\EKn$.
      Equivalently,  $J \subseteq \EKn$ is in the $G_0$-clustering if
and only if
there is a sequence $J_1,\ldots, J_i$ of subsets of $\EKn$ such that
each $J_j$ is
an elementary $G_0$-cluster,
$\bigcup_{j=1}^{i} J_j = J$, and $J_k\cap\bigcup_{i=1}^{k-1} J_j\ne \emptyset$
for
$2\le k
\le i$. (This definition of clusters corrects an error in the
definition in~\cite{W}.
The usage of it in~\cite{W} is consistent with the present definition.)

More generally, suppose $\Rset$ is any fixed set of nonempty graphs,
and information is
desired on the joint distribution of the subgraph counts for the
graphs in $\Rset$.
Then the appropriate clustering to consider is the minimal clustering
containing every elementary $G$-cluster for every $G \in \Rset$.  We call this the {\em
clustering generated by $\Rset$}. Of course, if $\Rset=\{G_0\}$, this is simply the $G_0$-clustering.

Henceforth in this paper we consider the clustering generated by a fixed
set of graphs $\Rset$, and assume that each graph in $\Rset$  has no isolated
vertices. Our first proposition considers a general set $\Rset$, and after that we restrict to only two kinds of clustering: the  $G_0$-clustering, and the one generated by
$\Rset=\{G_0,K_2\}$, which we call the {\em
$G_0^*$-clustering}. Note that a 1-element subset of $\EKn$ cannot
have a nontrivial proper intersection with any other cluster. It follows
that the $G_0^*$-clustering consists of the clusters of the
$G_0$-clustering, together with all the 1-element subsets of $\EKn$. We
assume in all cases that $|E(G_0)|\ge 2$.

      For $H \subs \EKn$, a {\em
cluster of $H$} is any cluster in  $\clusters$ contained in $H$.
A {\em maximal cluster} $Q$ of $H$ is cluster of $H$  which is contained in no
larger cluster of $H$.  Equivalently, $Q$ is a subset of $H$
such that
$Q\in\clusters$ and such that for every  $J\in\clusters$ with
$J\subseteq H$, either
      $J\subseteq Q$ or  $J\cap Q = \emptyset$.  (The case of nonempty
intersection is
excluded by the definition of a clustering.) For example, if $\clusters$ is the
$G_0$-clustering
      and $H$ is an arbitrary subset of $\EKn)$,
a maximal
cluster of
$H$ whose cardinality is
$|E(G_0)|$  must be an elementary $G_0$-cluster contained in $H$ having empty intersection with every other elementary $G_0$-cluster in $H$.  

Being a subset of $\EKn$, a cluster induces a subgraph of $K_n$. The
isomorphism class of the subgraph is called the {\em type} of the
cluster and also
of the subgraph. The set of types will be denoted $\types$, and we use $\type$
to denote the function which maps a cluster or the corresponding graph to its
type. Given
$t\in\types$, we use the notation
$|t|:=|\{S\subseteq
\EKn:\type(S)=t\}|$. Note that this depends on $n$, whereas $t$ is fixed. 

We will define a special
nonempty finite set $\smalls$ of types which is closed under taking subsets,
i.e.\ which satisfies
$$\mbox{if $S$, $S' \in \clusters$,  $\type(S) \in \smalls$ and $S' \subs S$
then
$\type(S') \in \smalls$}.$$
Let $s=|\smalls|$ be the number of types in $\smalls$.

The types in $\smalls$ will be called {\em small}, and any cluster $Q$ with
$\type(Q)\in
\smalls$ is also called small. Any type or cluster
which is not small is called {\em large}.
An {\em unavoidable} cluster is any large cluster which is a union of a small
cluster $Q$ and a set of small clusters all pairwise disjoint and all having
nonempty intersection with $Q$. The set of  types
of unavoidable clusters is denoted by $\unavoid$. (The term
``unavoidable" refers to the fact that large clusters created in a certain way, to be specified later, cannot avoid being in $\unavoid$.)

We will need to record how many subgraphs of every small type are present
in a given graph. So we consider the set $\functions$ of all
non-negative integer functions defined on
$\smalls$. For any $H \subs \EKn$, define $s_H$ to be the function in
$\functions$ such that,  for all $t \in
\smalls$,
$s_H(t)$ is the number of maximal clusters of $H$ of type $t$.
The function $\d_t\in\functions$ has value 1 at $t$ and 0 elsewhere.

All our basic work is in $\Gnp$, the standard edge-independent (binomial)
      model for
random graphs, and
$\pr$  and
$\ex$ denote probability and expectation in this space. $G$ denotes a random
graph in
$\Gnp$ and $q$ always denotes $1-p$.
      For $H \subs \EKn$, the event $H \subs E(G)$
is denoted by
$A_H$, so that
$\pr(A_H) = p^{|H|}$. The main objects we work with are, for each $f \in
\functions$, the set $\class_f$  consisting of graphs $G$ on $n$ vertices
containing no large clusters and such that $s_{E(G)} = f$. For $f \notin
\functions$, for example if
$f$ has a negative value on $\smalls$, we define
$\class_f=\emptyset$.
We write $\BP(\class_f)$ for $\BP(\Gnp\in\class_f)$.


For types $u,t\in \smalls$ and for $h \in \functions$, define, for
any fixed cluster $J$ of type $u$,
\eq\lab{cdef}
c(u,t,h) = \sum_{{Q\in \clusters \atop Q\subseteq J}\atop \type(Q)=t}
\sum_{H\cup Q=J\atop{s_H=h}} p^{|Q
\cap H|}q^{|J\setminus H|}.
\en
      Since the clustering generated by any set $\Rset$  is
symmetrical, $c(u,t,h)$ is clearly independent of the choice of $J$ with $\type(J)=u$. Note that in the special case $u=t$,
\bel{ctt}
c(t,t,h) =
\sum_{H\subseteq J\atop {s_H=h}} p^{|H|}q^{|J\setminus H|},
\ee
and in particular
\bel{ctt0}
c(t,t,\bz) = 1+O(p).
\ee

We use $\n(G)$ and $\m(G)$ for the numbers of vertices and edges of a graph $G$
respectively, and extend the notation to arbitrary subsets $H$ of $\EKn$,
so that $\n(H)$ is the number of vertices of the graph induced by $H$ and
$\m(H)$ is the number of edges. In particular,
this applies to clusters $H$.
We also use $\n(t)$ for the number of vertices in each cluster of type
$t$ and $\m(t)$ for the number of edges.

Let $[n]_k$ denote $n(n-1)\cdots(n-k+1)$.
For $t \in \types$, let $Q$ be any cluster of type $t$ and $|\mbox{aut}(Q)|$
      the number of automorphisms
of the graph induced by $Q$. Then
\bel{tsize}
|t|=\frac {[n]_{\n(Q)}}{|\mbox{aut}(Q)|},
\ee
and
\bel{lambdadef}
\mean_t := |t| p^{\mu(Q)}  =\Theta(n^{\nu(Q)}p^{\mu(Q)})
\ee
is the expected number
of different copies, in $G \in \Gnp$, of the subgraph induced by $Q$.

Our first result is obtained by simple counting. 
\begin{prop}\lab{basic}
For $f \in \functions$ and $t \in \smalls$,
$$
\frac{\pr(\class_{f+\delta_t})}{\pr(\class_{f})}=
\frac{\mean_t}{(f(t)+1)c( t,t,\bz)}
\left(
1-\Sigma- \frac{\theta(f,\delta_t)}{|t|\pr(\class_f)}
                                         \right)
$$
where
\bel{Sigmadef0}
\Sigma = \sum_{{u \in \smalls\atop h \in \functions}\atop (u,h)\ne (t,\bz)}
\frac{(f(u)-h(u)+1)c(u,t,h)\pr(\class_{f-h+\delta_u})}{\mean_t \pr(\class_f)}
\ee
and
\bel{first}
0\le\theta(f,\delta_t)\le \sum_{L\stcolon \type(L)\in\unavoid}\sum_{Q,H\subs
L\atop{\type(Q)=t
\atop L\setminus Q\subs H}} \pr(\class_{f-s_H})\left(\frac{p}{q}\right)^{|H|}.
\ee
\end{prop}
\proof
Note that
$$
c(u,t,h)p^{-\m(t)} = \sum_{{Q\in \clusters \atop Q\subseteq J}\atop
\type(Q)=t}\sum_{H\in\class_h\atop{H\cup Q=J}}
\frac{q^{|J\setminus H|}}{p^{|Q\setminus H|}}=
\sum_{{Q\in \clusters \atop Q\subseteq J}\atop \type(Q)=t}
\sum_{H\in\class_h\atop{H\cup Q=J}}
\left(\frac{q}{p}\right)^{|J\setminus H|},
$$
where we have used the fact that $J\setminus H=Q\setminus H$ follows
from $H\cup Q=J$.
Consider every pair $(E,Q)$ where $E$ is the edge set of a graph in
$\class_f$ and $Q$ is a cluster of type $t$. Classifying $E\cup Q$ according to
the type of its maximal cluster $L$ containing $Q$, and, in the case that
$\type(L)=u\in
\smalls$, subclassifying according to $h = s_{E\cap L}$, gives
\bel{second}
|t|\pr(\class_{f})=\Bigl(
\sum_{u \in \smalls\atop h \in \functions}
(f(u)-h(u)+1)c(u,t,h)p^{-\m(t)}\pr(\class_{f-h+\delta_u})\Bigr) +
\theta(f,\delta_t),
\ee
where the $\theta$ term is bounded as in the statement of the proposition. This
term comes from observing that if $L$ is a large cluster, then
it must be unavoidable since $E$ has no large clusters, and from
considering the
subset of
$\EKn$ obtained by removing the set $H$ of all edges of $E$ in $L$.
Multiplying~\eqn{second} by $p^{\m(t)}$ gives
$$
\mean_t\pr(\class_{f})=
\Bigl( \sum_{u \in \smalls\atop h \in \functions}
(f(u)-h(u)+1)c(u,t,h)\pr(\class_{f-h+\delta_u})\Bigr) +
\theta(f,\delta_t)p^{\m(t)}
$$
and rearranging the terms, isolating the one with $(u,h)=(t,\bz)$, finishes the
proof.
\qed

We now lay the groundwork for asymptotic results. Henceforth, we consider only the $G_0$- and $G_0^*$-clusterings for some fixed graph $G_0$ with at least two edges. Denote the set of proper subgraphs of $G_0$ which contain at least one edge by
$\extensions$. Recalling that $|E(G_0)| \ge 2$, we
    define the {\em extension value} of $G_0$ to be
\eq\lab{xdef}
x=x(G_0,p,n)=\max_{G_1\in\extensions}\, n^{\n(G_0)-\n(G_1)}p^{\m(G_0)-\m(G_1)}.
\en
For example, if $G_0$ is a triangle,
\bel{x}
x=\max(np^2, p,p^2)=\max(np^2, p).
\ee
      The significance of the extension value lies in the
fact that 
$n^{\n(G_0)-\n(G_1)}p^{\m(G_0)-\m(G_1)}$ is the asymptotically
important part of
$$
{\left({n-\n(G_1)}\atop {\n(G_0)-\n(G_1)}\right)} p^{\m(G_0)-\m(G_1)}.
$$
To interpret this quantity, first distinguish one of the subgraphs of $G_0$ isomorphic to $G_1$. For $G_2$ isomorphic to  $G_1$, conditional upon $G_2\subseteq \Gnp$, the quantity above is the expected number of isomorphisms from $G_0$ to a subgraph of $\Gnp$ that map the distinguished copy of $G_1$ onto $G_2$.

 For $H\subseteq E(K_n)$ define $\Phi(H,G_0)$  
to be the
expected number of subgraphs of $G\in\Gnp$ that are isomorphic to
$G_0$ and whose edge set contains $H$, conditional
      on $H\subseteq E(G)$. Given a nonempty $H\subseteq \EKn$ which
induces a proper
subgraph of $G_0$, it follows from the remarks above that $\Phi(H,G_0)$ is $O(x)$, since there is a bounded number of ways to distinguish one of the subgraphs of $G_0$ isomorphic to $G_1$.

Put a partial ordering on the set of types
by defining $t$ to be strictly less than $u$ in the poset,
denoted by
$t\prec u$, if, and only if, any cluster of
type $u$ properly contains a cluster of type $t$.
If $t\prec u$, then a cluster of type $u$ can be obtained from a cluster $Q$
of type $t$ by a finite sequence of non-disjoint unions with clusters
$Q_0, \ldots
,Q_k$
such that each
$Q_i$ is the edge set of a graph isomorphic to some
$G_i\in \Rset$ and $Q_i\not\subs Q\cup(\bigcup_{j=0}^{i-1} Q_j)$. (Note that, in the $G_0^*$-clustering, it must be that $G_i=G_0$ for all $i$.) Thus, for
$G\in \Gnp$ the expected number of clusters of type $u$ in
$E(G)$ can be bounded above by a finite sum whose terms are all of the form
$\mean_t\prod_{i=0}^k\Phi(H_i,G_i)$ where $H_i$ corresponds to the intersection of
$Q_i$ with $Q\cup(\bigcup_{j=0}^{i-1} Q_j)$.  Hence, from the
conclusion of the previous paragraph, provided $x=o(1)$ we have
\be\lab{meanratio}
\mbox{if $t\prec u$ then }\frac{\mean_u}{\mean_t}=O(x).
\ee

Henceforth in this paper, we assume that $G_0$ is strictly balanced, with at least two edges. Let $X$
be the number of copies of $G_0$
in the random graph $\Gnp$.
It follows easily from the definition  (\ref{xdef}) of $x$ that the constant
$\chi$
defined in~\eqn{chidef} is the smallest number such that $p=o(n^{-\chi})$ implies $x=o(1)$.
 Hence,   there are functions
$p=p(n)$ such that $\mean_{\type(G_0)}\to\infty$ while $x(G_0,p,n)=o(1)$.
We also assume henceforth that 
$p=p(n)$ is restricted so that
for some fixed $\kappa>\chi$,
\be\lab{pkappa}
p=n^{-\kappa+o(1)}.
\ee
This will be enough for our purposes in view of  Lemma~\ref{panything}.

Fix $\eps>0$ and let $ \kappa\ge \chi+\eps$. 
Since $\mu(G_1)<\mu(G_0)$ for all $G_1\in\extensions$,  the expression maximised in~\eqn{xdef} is at most $(n^\chi p)^{\m(G_0)-\m(G_1)}\le n^\chi p$. Thus, 
\be \lab{xbound}
x(G_0,p,n) = O(n^{-\eps+o(1)}).
\ee
See \cite{JLR} for a general introduction to  the considerations relevant here.  Note that 
\bel{plessx}
p\leq x
\ee
by definition, as shown by setting the graph $G_1$
in~\eqn{xdef} equal to $G_0$ minus an edge.

For our asymptotic results, we work with a particular set of small cluster types defined as follows:
\bel{smallsdef}
 \smalls =\{t: \nu(t)/\mu(t)\geq\kappa \}.
 \ee
Then for
$t\in\smalls$, the expected number $\lambda_t$ of subgraphs
of type $t$ is bounded below by $\lambda_t\geq n^{-o(1)}$ (here the negative sign is not necessary, just indicative, since $o()$ bounds the absolute value), 
since by~\eqn{tsize},~\eqn{lambdadef} and~\eqn{pkappa}, 
\bel{lasize}
\la_t=\Theta(n^{\nu(t)-\kappa \mu(t)+o(1)}).
\ee
The set $\smalls$
is finite by (\ref{xbound}) and (\ref{meanratio}).
Hence, defining
\bel{maxlarge}
\mean_\larges:=\sup_{t\notin \smalls}\mean_t
\ee
we obtain
\bel{maxlargebound}
\mean_\larges=O(n^{-\epsilon'})
\ee
for some $\epsilon^\prime>0$
by our definition of $\smalls$.
While we are at it, due to a technicality we assume $\kappa<2$, so that  
$p$  satisfies the very weak growth condition
      \bel{Egrows}
      n^2p>n^{\epsilon^{\prime\prime}}
      \ee
for some $\epsilon^{\prime\prime}>0$.\ This ensures that the number of edges in the random graph tends to
infinity at a reasonable rate.
 Imposing this condition is without loss of generality, since  the omitted case follows from the case considered. For example, 
 the $p$ such  that $p\sim n^{-c\nu(G_0)/\mu(G_0)}$  are 
covered for all 
 $1<c<2\mu(G_0)/\nu(G_0)$, and this is well below the threshold of appearance of copies of $G_0$. Hence, each term in the power series must tend to zero for such $c$, and must also tend to 0 when  $\kappa\ge  2$. 
 The assumption $\kappa<2$ also ensures that, in the case of the $G_0^*$-clustering, the single edge cluster is in $\smalls$.  Note that if $n^2p=o(\sqrt n)$,  the random graph is in any case not interesting, as it is asymptotically almost surely a matching. 

Define
$$
\smalls_0 =\{t\,:\, \nu(t)/\mu(t)=\kappa\},\quad \smalls_1 =\smalls
\setminus \smalls_0,
$$
\bel{mtdef}
m_t = \left\{
\begin{array}{ll}
3\lambda_t&\mbox{ if }t\in\smalls_1\\
\lambda_t\log n&\mbox{ if }t\in\smalls_0.
\end{array}
 \right.
\ee
Note that  $\smalls_0$ will often be empty, but if it is nonempty, the types in $\smalls_0$ are the rarest types of small clusters in the random graph, and for  $t\in \smalls_0$, we have  $\la_t=n^{o(1)}$  and  hence $m_t=n^{o(1)}$.  Any type in $\smalls_0$ is   maximal in $\smalls$ by~\eqn{meanratio}.  Thus, for later reference we may note that,
for some positive $\eps'''$, 
\bel{s0s1}
\la_t>n^{\epsilon'''} \mbox{ for } t\in \smalls_1, \quad
\la_t = n^{o(1)}  \mbox{ for } t\in \smalls_0.
\ee

Let $\smallfunctions=\smallfunctions(n)$
be the set containing those functions $f\in\functions$ such that for all $t
\in \smalls$,
\eq\lab{fbound}
f(t)\le
m_t.
\en
For integer-valued $h$ with $f,f+h\in\functions$, we define
\bel{rhodef}
\rho(f,h)=\frac{\pr(\class_{f+h})}{\pr(\class_f)}
\ee
and  for $t\in\types$, $f\in\functions$ define
\bel{gammadef}
\ga(f,t)=\frac{\rho(f,\delta_t)(f(t)+1)}{\lambda_t}.
\ee
The motivation for focussing on $\ga$ is that if the numbers of clusters of the
various small types were independent Poisson variables, then all the $\ga$'s
would be exactly 1. Proving that they are close to 1 shows that the variables
are approximately Poisson. We will be measuring the difference between
the Poisson probability and the true probability of $\class_f$  very
accurately for some values of $f$.

 Ultimately, we wish to estimate $\ga(f,t)$, and will achieve this in Corollary~\ref{expanded}. The proof is complicated, so is broken up into several parts, obtaining progressively simpler approximations. The downside of breaking it up like this is that it requires repeating the same kinds of inductive arguments several times. We first  obtain a more useful bound on the function $\theta(f,\d_t)$  appearing in Proposition~\ref{basic}. Let $t^*$ denote the type of the single edge cluster, which of course only appears in the $G_0^*$-clustering.

\begin{prop}\lab{thetabound2}
Uniformly for every $f\in\smallfunctions$
and every
$t\in \smalls $,
$$
\frac{\theta(f,\delta_t)}{|t|\BP(\class_f)}
= O\left(\frac{ \phi_t\mean_\larges}{\mean_t}\right), 
$$
where $ \phi_{t}= n^{o(1)}$ for $t=t^*$ and  $\phi_{t}=1$ otherwise.
 Moreover, for all $f\in\smallfunctions$
and $t\in \smalls$, uniformly,
$$
\ga(f,t)=1+O(xn^{o(1)}).
$$
\end{prop}
\noindent
{\bf Note.\ }  The proof will reveal that the factor $n^{o(1)}$ can be replaced by  the maximum of $f(t')/\mean_{t'}$ for $t'\in\smalls$, which is always at most $\log n$. However, $n^{o(1)}$ is tight enough for our purposes here. 
Also, $\mean_\larges$ can be replaced by the maximum value of  $\mean_u$ over all $u\in\larges$ such that  $t\prec u$.
\smallskip

\proof 
In this   proof, as in the proposition's statement,  the constants implicit in the $O()$ terms depend only on the choice of clustering and $\kappa$, as do the bounds implicit in the notation $\sim$ and $o(1)$. We will use induction on
$f\in \smallfunctions$. Order $\smallfunctions$ lexicographically; that is $g<f$ if, and only if,
$g\neq f$ and $g$ has a smaller value than $f$ in the first component at
which they differ. This induction is crucual to the whole approach of this paper, and is rather unusually complex, since  for the $G_0^*$-clustering,  the induction actually begins with the graph on $n$ vertices and no edges.  So we formulate a statement that pays explicit attention to the implicit constants in $O()$: what
we claim is that there exists   constants $C$ and $C'$, a number $N_0$ and a function $1\le \phi^* = \phi^* (n)=n^{o(1)}$
 (all depending only on the clustering and $\kappa$) such that, for $n\ge N_0$ and all relevant $f$ and $t$,
\bel{Cbound}
\frac{\theta(f,\delta_t)}{|t|\BP(\class_f)}
\le C  \phi_t\,  \frac{ \mean_\larges}{\mean_t}, 
\ee
where $\phi_t=\phi^*$ for $t=t^*$ and $\phi_t=1$ otherwise, and furthermore
\bel{Cprimebound}
|\ga(f,t)-1|\le C' \phi_t x\le 1/2.
\ee
To prove this, we can assume that for this particular $C$, and $n$ large enough, these inequalities hold when $f$ is replaced by any $g<f$ (in the lexicographic ordering).

We first discuss the bound 
involving $\theta$. Here,  by~\eqn{maxlarge}, it is enough to show the bound 
$C\phi_t\mean_{\type(L)}/\mean_t $ where  $\type(L)\notin\smalls$ (which then justifies the second part of the note after the statement of the proposition). Moreover, of~\eqn{Cprimebound} we will only use the inequality
\bel{Cprimebound2}
|\ga(f,t)-1|\le   1/2.
\ee
Since the number of clusters of the complete graph $K_n$ which are
isomorphic to a given $L$ is
$O(n^{\nu(L)})$, and since the number of types of unvoidable clusters
is by definition bounded, we may use~\eqn{first} and $q\sim 1$ to
obtain the bound  
\bel{thbound}
\frac{\theta(f,t)}{|t|\pr(\class_f)}=
O(1)\max_{ {{\type(L)\in \unavoid}\atop{\type(Q)=t, Q\subseteq L}} \atop
{L\setminus Q \subseteq H \subseteq L} }
n^{\n(L)-\n(Q)}p^{|H|} \frac{\pr(\class_{f-s_H})}{\pr(\class_f)} 
\ee
for $n$ sufficiently large (which in particular ensures that  $\pr(\class_{f-s_H})\ne 0$).  Here, recalling~\eqn{lambdadef} we see that
\bel{facts}
|H|\ge |L|-|Q|, \quad \mean_t=O(n^{\n(Q)}p^{|Q|}), \quad
n^{\n(L)}p^{|L|}=O(\mean_{\type(L)}).
\ee

In the case $f=\bz$,   we may assume $s_H=\bz$ in~\eqn{thbound}, since otherwise, $\class_{f-s_H}$ is empty. Thus, by~\eqn{facts},
we have the bound $O(\mean_{\type(L)}/\mean_t)$ on each term in~\eqn{thbound}. Since $\type(L)\notin\smalls$, we are done in this case.

In the case $\bz \ne f\in\smallfunctions$, suppose the claim has been
shown when $f$ is replaced by any $g<f$. We need to show that, when $C$ is large enough, the very same $C$ applies in the statement for $f$.  Denoting a general term in the maximum in~\eqn{thbound} by $M$, since $\type(L)\in \unavoid \subseteq \larges$,  it suffices to show   that  $M=O(\mean_{\type(L)}/\mean_t)$, or $M=O(n^{o(1)}\mean_{\type(L)}/\mean_t)$ in the case of the $G_0^*$-clustering  (and then choosing $\phi^*$ appropriately). We may write
\bel{rhoprod}
\frac{\pr(\class_{f-s_H})}{\pr(\class_f)}=
\prod_{i=1}^{k} \rho(f_i\, ,-\delta_{u_i})
\ee
for some sequence
$u_1,u_2,\ldots , u_k$ in $\smalls$
such that $\sum_{i=1}^k\delta_{u_i}=s_H$ and  where
$f_i=f - \sum_{j=1}^{i-1}\delta_{u_j}$.
By definition, an unavoidable cluster has size at most $r(r-1)$ where 
$r$ is the size of the largest small cluster. Hence, the upper index $k$ in the above product is at most $r(r-1)$.  Note also that each $f_i$ occurs before $f$ in the lexicographic order, and~\eqn{Cprimebound2} inductively implies $1/2\le \ga(f_j-\d_{u_j},\d_{u_j})\le 3/2$
for all $j\ge1$.
Note that
$$
\rho(f_j,-\delta_{u_j}) =
\frac{1}{\rho(f_j-\delta_{u_j},\delta_{u_j})}
=
\frac{f_j(u_j)}{\lambda_{u_j}\ga(f_j-\delta_{u_j},u_j)}.
$$
Suppose firstly that, in~\eqn{rhoprod},  $u_i\in\smalls_1$  for all $i$. Then by~\eqn{fbound},
$f_j(u_i)/\lambda_{u_i} \leq 3$ for all $i$, and by~\eqn{Cprimebound2}  inductively
$\ga(f_i-\delta_{u_i},u_i)^{-1}\le 2$, so we deduce that the product in~\eqn{rhoprod} is   $O(1)$.   Now~\eqn{facts} implies that   
$M=O\left(\mean_{\type(L)}/\mean_t\right)$, 
as required.

Suppose on the other hand that, for  some term
in~\eqn{rhoprod}, 
there is some $j^\prime$ for which
$  u_{j'}\in \smalls_0$.
Recall that $\mean_{ u_{j'}}=n^{o(1)}$ by~\eqn{s0s1},  and hence
\be\lab{rhoprodbound}
\frac{\pr(\class_{f-s_H})}{\pr(\class_f)}
      =\rho(f,-s_H)
      =O(n^{o(1)})
\ee
using the same argument as for analysing~\eqn{rhoprod} above.
Also note that
\be\lab{otherstuff}
n^{\n(L)-\n(Q)}p^{|H|} =
n^{\n(L)-\n(Q)} p^{|L\setminus Q|} p^{|H\cap
Q|}=O(\lambda_{\tau(L)}/\lambda_t)p^{|H\cap Q|}.
\ee
There are two subcases to consider.
Firstly, if $|H\cap Q|\geq 1$, then
$p^{|H\cap Q|}n^{o(1)}\le pn^{o(1)} =o(1)$ and hence 
$M=O(\lambda_{\type(L)}/\mean_t)$ as required.
The second subcase is $|H\cap Q|=0$.
Then $H$ contains a cluster $Q'$ of type $u_{j^\prime}$, disjoint from $Q$.
It follows that there is a sequence $Q_1,\ldots, Q_\ell$ of elementary clusters, each nontrivially intersecting the next, with $Q_1\cap Q'\ne \emptyset$ and $Q_\ell\cap Q\ne \emptyset$, $Q_\ell\ne Q$. 
 We will consider two subsubcases  of this second case.

Suppose firstly  that $Q\not\subseteq Q_\ell$, and so $Q'':= Q'\cup\bigcup_{i=1}^\ell Q_i$ is a cluster satisfying $Q'\subset Q''\subset L$,
where the inclusions are proper and $\tau(Q')=u_{j^\prime}$.
It follows by (\ref{meanratio}) and (\ref{s0s1})
that $\lambda_{\tau(Q'')}=O(\lambda_{u_{j^\prime}}x) =O( n^{o(1)}x)$ since $u_{j'}\in \smalls_0$. Thus $\tau(Q'')\in\larges$, and hence  by   the definition~\eqn{maxlarge} of $\la_\larges$, we have
$\lambda_{\tau(Q'')}\le \la_\larges$. Similarly,   $\lambda_{\tau(L)}=O( x\lambda_{\tau(Q'')} )= O( x\lambda_\larges )$, and now using~\eqn{rhoprodbound} and~\eqn{otherstuff}   in~\eqn{thbound}  gives $M=O\left(x\mean_{\type(L)} n^{o(1)}/\mean_t\right)  =O\left(\mean_{\type(L)} /\mean_t\right)$ as required.

For the other subsubcase  $Q\subseteq Q_\ell$, recall that $ Q_\ell\ne Q$. As $Q_\ell$ is elementary, it follows that this can only occur for the  $G_0^*$-clustering, and $Q$ must be a single edge (and  its type $t$ equals $t^*$).  Using~\eqn{rhoprodbound} and~\eqn{otherstuff}   in~\eqn{thbound}  gives $M   =O\left(\mean_{\type(L)} n^{o(1)}/\mean_t\right)$ in this case, as required. We note that in fact the bound can be strengthened to $O\left(\mean_{\type(L)}/\mean_t\right)$  unless $Q_\ell=L$, $\ell=1$ and $j=1$, and looking back at the above argument, we may use $ f_j(u_j)/ \lambda_{u_j}$ in place of $n^{o(1)}$, as noted after the proposition's statement.

 We turn now to proving the bounds 
$$
|\ga(f,t)-1|\le C'{\phi_t} x 
$$
for all $t\in\smalls$, and here we may assume by induction that~\eqn{Cprimebound2} holds with $f$ replaced by any $g<f$, and that, as we have just shown,~\eqn{Cbound} holds.
We also know that $c(t,t,\bz) = 1 + O(p)$ from \eqn{ctt0}.
So it suffices to
show that $\Sigma$ in the statement of the Proposition \ref{basic} is
$O(\phi_t x)$.
 Since
$\smalls$ is fixed, there is a bounded number of
terms in the sum, and each may be written as 
\bel{term}
\ga(f-h,u)\,\frac{\lambda_u}{\lambda_t}\,
c(u,t,h) \rho(f,-h).
\ee 
Note that the  argument that produced~\eqn{rhoprodbound} gives, in this case, 
$\rho(f,-h) = O(n^{o(1)})$. So  (again by appropriate choice of ${\phi^*}$) we only need to show that the product of the remaining factors in~\eqn{term} is $O(xn^{o(1)})$.
 
Let $\hfunctions$ denote the set of $h\in\smallfunctions$ for which
there are $t,u\in\smalls$ such that $c(u,t,h)\neq 0$.
Note that the cardinality of $\hfunctions$ is bounded.

Inside the present main inductive step, we use a second level of induction on $t$, going from greatest to smallest in the relation `$\prec$'. Assume first that $t$ is maximal. Since $u\in\smalls$, it is necessary that $u=t$ and $ h\ne \bz$ for such a term to be included in $\Sigma$. Then   $\ga(f-h,t)\le 3/2$ by~\eqn{Cprimebound2} inductively.
 Furthermore, since the graphs in $\Rset$ are nonempty and
$H\ne\emptyset$ in~\eqn{ctt}, we have
$c(t,t,h)=O(p)=O(x)$, which gives the
desired result.

Suppose next that $t$ is not maximal.
A term~\eqn{term} with $u=t$ and $h\neq \bz$
is $O(xn^{o(1)})$ for reasons as in the  previous  paragraph.    On the other hand, 
for $u\neq t$ and $h\in\hfunctions$, clearly $c(u,t,h)=O(1)$.
If $c(u,t,h)\ne 0$, then by the definition~\eqn{cdef}, $t\prec u$, and then  
$\gamma(f-h,u)\le 3/2$ by~\eqn{Cprimebound2}  inductively, and $\lambda_u/\lambda_t=O(x)$ by~\eqn{meanratio}. Once again,~\eqn{term} is  $O(xn^{o(1)})$. For appropriate choice of ${\phi^*}$ and $C'$, we now have $|\ga(f,t)-1|\le C'{\phi_t} x $. Thus, in view of the bound~\eqn{xbound} on $x$, for appropriate choice of $N_0$, we have~\eqn{Cprimebound} in full. This completes the inductive step, and ~\eqn{Cbound} and~\eqn{Cprimebound} imply the lemma.
\qed
It is useful to rewrite Proposition \ref{basic}
in terms of the $\gamma$'s.
It says that for $f \in \functions$ and $t \in \smalls$,
\bel{gammaalt}
\ga(f,t)=\frac{1}{
c(t,t,\bz)}\left(
1-\Sigma
- \frac{\theta(f,\delta_t)}{|t|\pr(\class_f)}
\right),
\ee
where $\Sigma$ is defined  by \eqn{Sigmadef0}. Writing 
$$\frac{\pr(\class_{f-h+\delta_u})}{ \pr(\class_f)}
=\frac{\pr(\class_{f-h})}{ \pr(\class_f)}\cdot
\frac{\pr(\class_{f-h+\delta_u})}{ \pr(\class_{f-h})}
$$
and using~\eqn{rhoprod} for the first factor gives
\bel{Sigmadef}
\Sigma = \sum_{{u \in \smalls\atop h,f-h  \in \functions}\atop
(u,h)\ne (t,\bz)}
\frac{\lambda_u}{\lambda_t}c(u,t,h)\ga(f-h,u)
\prod_{i=1}^{k} \frac{f_i(t_i)+1}{\lambda_{t_i}\ga (f_i,t_i)},
\ee
which is  a function of $f$ and $t$, where, for each $h$, $t_i$, $i=1,\ldots,k$ is a sequence in $\smalls$ such that
$h=\sum_{i=1}^{k}\delta_{t_i}$ and
$f_i=f-\sum_{j=1}^{i}\delta_{t_j}$. Here and henceforth, we may
choose a canonical sequence $t_1, \ldots , t_k$ for each $h$ such
that $c(u,t,h)\ne 0$ for some $u,t\in \smalls$. Note that $k$ is bounded
because
$\smalls$ is finite.

Approximations to the $\ga$'s may be defined recursively by ignoring
the term  containing   $\theta (f,\delta_t)$ in  (\ref{gammaalt}). Thus,   we define:
$$
      \hamma(f,t)=
\frac{1}{
c(t,t,\bz)}\left(
1-\hat\Sigma
\right)
$$
where
\bel{hammadef}
\hat\Sigma = \sum_{{u \in \smalls\atop h,f-h \in \functions}\atop
(u,h)\ne (t,\bz)}
\frac{\lambda_u}{\lambda_t}c(u,t,h)\hamma(f-h,u)
\prod_{i=1}^{k} \frac{f_i(t_i)+1}{\lambda_{t_i}\hamma(f_i,t_i)}
\ee
is a function of $f$ and $t$.
\begin{prop}\lab{hammaclose}
Uniformly for all $f\in\smallfunctions$ and $t\in\smalls \setminus\{t^*\}$,
$$
|\hamma(f,t)-\ga(f,t)|= O\left(\frac{x+\phi_t\lambda_\larges}{\lambda_t}\right)
$$
 where $ \phi_{t}= n^{o(1)}$ for $t=t^*$ and  $\phi_{t}=1$ otherwise.
\end{prop}
\proof We use an inductive scheme as we did for Proposition~\ref{thetabound2}. The  initial step of the
outer induction is $f=\bz$, and the initial step  of the inner induction has $t$ maximal in $\smalls$. The initial steps are considered below.  

 We aim
to show inductively that
\bel{inductgamma}
\ga(f,t)=\hamma(f,t)+ O_t\left(\frac{x+\phi_t\lambda_\larges}{\lambda_t}\right).
\ee
where $O_t()$ denotes $O()$ with the implicit constant depending on $t$.
(Although this implies the same statement for a uniformly defined implicit
constant, the induction argument requires different constants for each $t$,
larger constants for ``smaller" $t$. Constraints on the sizes of these constants are
implicitly determined in the proof below.) 
By~\eqn{ctt0}, the definition~\eqn{x} of $x$, and Proposition~\ref{thetabound2}, it suffices
to show 
\bel{inductsigma}
\Sigma=\hat\Sigma+ O_t\left(\frac{x+\phi_{t}\lambda_\larges}{\lambda_t}\right).
\ee

Instead of proceeding step by step through the induction, the argument is
made by focussing on the relevant considerations for  an arbitrary step, whether it be an initial
step (for $f$ or for $t$) or an arbitrary inductive step.

First, notice that if some $t_i=u$  in (\ref{hammadef}), then it must
be that $k=1$, $h=\delta_u$, $f_1=f-h$ and the $\hamma$'s cancel.
This means that the corresponding terms in
$\Sigma$ and $\hat\Sigma$ are equal, so henceforth
whenever $k\ge 1$, we may assume that $t_j\prec u$ for all $j$.

If $h=\bz$ in a term in $\Sigma$, or $\hat
\Sigma$, then the value of $k$ in that term is 0, and the product in that term is empty, and equal to 1. On the other
hand, suppose that
$h\ne \bz$. As shown above, we may assume  that each $t_j\prec u$.
Thus, in (\ref{hammadef}), $\lambda_{t_j}\to\infty$ for all $j$,
because
if any of these were bounded, it would imply $\mean_u =
O(x)$ and so $u\notin\smalls$. The ratios
$(f_i(t_i)+1)/\lambda_{t_i}$ in~\eqn{Sigmadef} and~\eqn{hammadef} are
therefore
$O(1)$ by~\eqn{fbound}. We have from Proposition~\ref{thetabound2} 
 that
$\ga(f,t)\sim 1$  uniformly, and it
is also immediate that $c(u,t,h)=O(1)$, and
$1/c(t,t,\bz)=O(1)$ by~\eqn{ctt0}. The combination of these facts 
shows that each
$\hamma(f_i,t_i)$ in~\eqn{hammadef} is $1+o(1)$, with the convergence uniform
over all $f_i$ and $t_i$. This implies in particular that the product
in~\eqn{hammadef} is in all cases
$O(1)$.

We will estimate the difference between the summands in~\eqn{Sigmadef}
and~\eqn{hammadef} using 
\bel{deltas}
(A+\delta_A)(B+\delta_B)-AB=O(|\delta_A B|  +
|A \delta_B|),
\ee
which holds provided that $\delta_A=O(A)$ or $\delta_B=O(B)$.
  We will show that for
$(u,h)$ as in the scope of the summation in~\eqn{Sigmadef},
\eq\lab{diff1}
\left|\hamma(f-h,u)-\ga(f-h,u)\right|\frac{\lambda_u}{\lambda_t}c(u,t,h)
=
   \left\{ \begin{array}{ll} \displaystyle \raisebox{-0.6cm}{\rule{0cm}{1cm}}
O_u\left(\frac{x+\lambda_\larges}{\lambda_t}\right)&\mbox{ if
}t\prec u\\ \displaystyle
O_t\left(x\frac{x+\phi_t\lambda_\larges}{\lambda_t}\right)&\mbox{ if
}u=t,\end{array} \right.
\en
and, for factors appearing in the product in~\eqn{hammadef}  with $\ t_i\prec
u$,
\eq\lab{diff2}
\left|\hamma(f_i,t_i)-\ga(f_i,t_i)\right|\frac{\lambda_u}{\lambda_t}
=O_{t_i}\left(x\cdot \frac{x+\phi_{t_i}\lambda_\larges}{\lambda_t}\right).
\en
In view of the above observations, these imply
$$\Sigma=\hat\Sigma+
\sum_{{u\in\smalls:\,t\prec u}}
\frac{x+\phi_{t}\lambda_\larges}{\lambda_t}O_u(1)
+\sum_{v\in\smalls}
\frac{x+\phi_{t^*}\lambda_\larges}{\lambda_t}O_v(x).
$$
Equation~\eqn{inductsigma}
will then follow, since the summations contain a bounded number of terms,
 and in the first summation the constant
implicit in $O_u()$ may be used in defining the constant implicit in
$O_t()$, whilst in the second summation the bound is 
 $o\big((x+\phi_t\lambda_\larges)/\lambda_t\big)$ by induction   using $x\phi_{t^*} =o(1)$). Note that for the initial step of the inner induction, when $t$ is
maximal in $\smalls$, it  must be that $u=t$.

For each term in~\eqn{Sigmadef} and~\eqn{hammadef} we have  $(u,h)\ne (t,\bz)$, so the inductive
statement~\eqn{inductgamma} implies
$$
\left|\ga(f-h,u)-\hamma(f-h,u)\right|\frac{\la_u}{\la_t}
=O_u\left(\frac{x+\phi_u\lambda_\larges}{\lambda_t}\right).
$$
Note   that   $t\prec u$ implies $u\ne t^*$ and hence
$\phi_u=1$. Recalling $c(u,t,h)=O(1)$, and noting that in particular $c(t,t,h)=O(x)$  when $t=u$   (as $h\ne
\bz$ in that case), we have~(\ref{diff1}).
By the outer induction (which is on $f$) using~\eqn{inductgamma}, the left side of (\ref{diff2}) is of
order
\eq\lab{tibound}
 O_{t_i}\left(\frac{x+ \phi_{t_i}\lambda_\larges}{\lambda_{t_i}}\,
\frac{\lambda_u}{\lambda_t}
\right)=O_{t_i}\left(\frac{x+ \phi_{t_i}\lambda_\larges}{\lambda_{t}}\,
\frac{\lambda_{u}}{\lambda_{t_i}}\right)
\en
and by (\ref{meanratio}) and (\ref{xbound}) (noting that $t_i\prec u$ as
discussed above),
$\lambda_u/\lambda_{t_i}=O(x)$, which
completes the proof.
\qed

A recursive calculation of $\hamma$ using its definition, including (\ref{hammadef}),  would need to keep track of $\hamma(f,t)$ for each
$f\in\smallfunctions$ and $t\in\smalls$. By making further approximations, we may obtain a simpler
recursion for functions which are explicitly defined in a compact form, and not depending on
$f$. Recalling that $|\smalls|=s$, without loss of generality we denote $\smalls$ by
$[s]=\{1,\ldots ,s\}$. (Thus $t\in \smalls$ is represented by an integer. We apologise to the reader for the possible confusion resulting; in particular the definition (\ref{tsize}) of the function $|t|$, where $t$ is a type,
overrides the notation for absolute value of the integer. It only appears once or twice more.) 
The simpler recursion will define  $\hatma_t\in  \real[[n,p,g_1,\ldots, g_s]]$, i.e.\ a
formal power
series in $n$, $p$ and $g_1,\ldots,g_s$ with real coefficients. Occasionally it will be useful to regard $\hatma_t$ also as an element of $\real[[n,p]][[\bfg]]$ where $\bfg=(g_1,\ldots,g_s)$, meaning a formal power series with indeterminates  $g_1,\ldots,g_s$ and coefficients in $\real[[n,p]]$.
Later,  we will calculate the new estimates of $\gamma(f,t)$ by    setting   $g_i=f(i)/\la_i$ in  $\hatma_t$ for each $i$.

Note that $c(u,t,h)$ is a polynomial in $p$, and $1/{c(t,t,\bz)} =1+O(p)$
and can be expanded as power series  in  $p$. Also, by (\ref{lambdadef}),
for $t\prec u$,
$\lambda_u/\lambda_t$ is a polynomial in $n$ and $p$
with terms of the form $p^{\mu(u)-\mu(t)}n^i$, and, since $\mu(u)>\mu(t)$,
$\lambda_u/\lambda_t$ has zero constant term.
    With these interpretations, we will define 
$\hatma_t=\hatma_t(n,p,\bfg)\in  \real[[n,p,g_1,\ldots, g_s]]$ using
\eq\lab{hatmadef}
\hatma_t= \frac{1}{c(t,t,\bz)}\left(
1-\sum_{{u \in \smalls\atop h \in \functions}\atop (u,h)\ne (t,\bz)}
\frac{\lambda_u}{\lambda_t}c(u,t,h)\hatma_u
\prod_{i=1}^{k} \frac{g_{t_i}}{\hatma_{t_i}}
\right), \quad \hatma_t(0,0,\bz)=1
\en
simultaneously for all $t\in\smalls$, where the $t_i$ are defined as in~\eqn{Sigmadef}. Since $c(t,t,h)=O(p)$ for
$h\ne \bz$ and $(\la_u/\la_t)c(u,t,h)$ has zero constant term for $u\neq
t$, there is a unique set of formal power series
$\hatma_t(n,p,\bfg)$, $t\in \smalls$, defined by~\eqn{hatmadef}, and they all have constant term 1. It
will also be useful to rewrite (\ref{hatmadef})   as
\bel{wrecurse}
\hatma_t=
1+w_0(t)-\sum_{{u \in \smalls\atop h \in \functions}\atop (u,h)\ne (t,\bz)}
w(u,t,h)\hatma_u
\prod_{i=1}^{k}\frac{1}{\hatma_{t_i}},
\ee
\bel{wdefs}
      w_0(t)=\frac{1}{c(t,t,\bz)}-1,\quad
      w(u,t,h)= \frac{\lambda_uc(u,t,h)}{\lambda_tc(t,t,\bz)}
\prod_{i=1}^k g_{t_i}.
 \ee
Here~\eqn{wrecurse}  defines $\hatma_t$ as a power series in the $w$'s, which, if
substituted appropriately as  power series in  $n$, $p$ and $\bfg$ using~\eqn{wdefs}, results
in the same series as defined in~\eqn{hatmadef}.

Given a function $f\in\functions$, with a slight abuse of notation, define
\bel{hatmaf}
\hatma_t(f)=\hatma_t(n,p,\tilbfg)
\ee
where
$$
\tilbfg =  (f(1)/\lambda_1,\ldots,f(s)/\lambda_s).
$$
Thus, given $n$ and $p$, $\hatma_t(\cdot)$ maps functions $f\in\functions$ to numbers, whereas $\hatma_t$ is a power series.

Returning to our original setting,   $f\in
\smallfunctions$ (as defined at~\eqn{fbound}), and  $p$ is a function of $n$
such  that
$x=x(n,p)=O(n^{-\epsilon})$ by~\eqn{xbound}.
It might help to observe at this point that, for given $n$, $p$ and
$f$ satisfying these constraints, there is a unique value of $\hatma_t(f)$
determined from the equations~\eqn{hatmadef} and~\eqn{hatmaf},
as long as $n$
   is large enough. One way to prove this is to consider an
initial approximation for each $\hatma_t(f)$, and then, iterating the
approximations using~\eqn{hatmadef}, with $g_t$ set equal
to $f(t)/\lambda_t$, the current values of $\hatma_t$ on the right side
giving rise to updated values on the left side. This determines
a contractive mapping on the vector whose entries are $\hatma_t(f)$
($t\in\smalls$) which has a fixed point near the initial approximate solution
determined by  $\hatma_t(f)=1$ for all $t$.  
To flesh this out, we first examine the   definition of $\hatma_t$ in order to bound the error of  approximations.  Recalling~\eqn{pkappa} and~\eqn{xbound}, we have the following
lemma.

First, given particular values of $n$, $p$ and $f$, we define
$$
\tilde g_{t}=f(t)/\la_t,
$$
so that $\tilbfg = (\tilde g_1, \ldots, \tilde g_s)$, and let  $\tilde w(u,t,h)$ denote the  value of    $w(u,t,h)$ obtained if we  replace $g_{t_i}$ by $\tilde g_{t_i}$ in~\eqn{wdefs}. For convenience, similarly set $\tilde w_0(t)= w_0(t)$.  Recall that
$p$ has been assigned a function  of $n$ satisfying~\eqn{pkappa}, which is significant when considering issues of uniformity.

\begin{lemma}\lab{doubleyous}
Suppose that
$0\le \tilde g_t=\tilde g_t(n)=O(n^{o(1)})$, with $\tilde g_t(n)=O(1)$ if
$t\in\smalls_1$. Then
$\tilde w_0(t) =O(p)$ and
$\tilde w(u,t,h) = O(x)$ for each term in~\eqn{wrecurse}, where the bounds in
the $O()$ terms are uniform.
\end{lemma}

\proof
   From~\eqn{ctt0}, $\tilde w_0(t) =O(p)$ and, recalling that $k$ is bounded
in~\eqn{wdefs} and that $c(t,t,\bz)\sim 1$,
\bel{wbound}
\tilde w(u,t,h)=O\left(\frac{\mean_uc(u,t,h)}{\mean_t}  (\max_i \tilde g_{t_i})^k
\right).
\ee

Firstly, if $h=\bz$, then $k=0$, and $u\succ t$ by the condition in the
summation. So
$\tilde w(u,t,h) = O(x)$ by~\eqn{meanratio}.

Secondly, suppose that $h\ne \bz$ and $u=t$.
If $h=\delta_{t^*}$ (recall that $t^*$ is the type of the single-edge cluster), then  $c(u,t,h)=p$. By~\eqn{Egrows}, we have $t^*\in \smalls_1$. So, using the hypothesis of this lemma,
the maximum in~\eqn{wbound} is
$O(1)$, and thus
$\tilde w(u,t,h)=O(p)=O(x)$. In all other cases, if $c(u,t,h)\ne 0$
then~\eqn{ctt} gives $c(t,t,h)=O(p^2)$ since $s_H=h$ implies $|H|\ge
2$. By~\eqn{wbound}, again $\tilde w(u,t,h)=O(x)$.

Lastly, suppose that $h\ne \bz$ and $u\succ t$.
Here $\mean_u/\mean_t=O(x)$ by~\eqn{meanratio}, and so
we are done if the maximum in~\eqn{wbound} is $O(1)$. But this must happen
unless $t_i\in\smalls_0$ for some $i$. Since $H$ contains only subclusters of a cluster
of type $u\in \smalls$,~\eqn{meanratio} shows that this requires $t_i=u$.
Then we have $h=\delta_u$, and hence in~\eqn{cdef}, $Q\subseteq J$ and  $|Q\cap H|\ge 1$, and so
$c(u,t,h)=O(p)=O(x)$. Since the maximum in~\eqn{wbound} is $O(n^{o(1)})$, the
bound obtained is $O(x^2n^{o(1)})$, and the result follows in this case
also.
   \qed

Recall that $\hatma_t(f)$ is a function of $n$, $p$ and $f$.
\begin{lemma}\lab{converge}
For
$f\in
\smallfunctions$
and $p$ satisfying~\eqn{pkappa}, the series definition of
$\hatma_t(f)$ in~\eqn{hatmaf} converges absolutely for $n$ sufficiently large, and
$\hatma_t(f)=1+O(x)$, where the bound in the $O()$ notation is uniform.
\end{lemma}
\proof
For any $t\in\smalls_0$, it follows from the definition of $\tilde g_t$,
the upper bounds~\eqn{mtdef} and \eqn{fbound} on $f(t)$, and the asymptotics~\eqn{pkappa} of $p$,
that
$\tilde g_t=O(n^{o(1)})$.
On the other hand, if $t\in\smalls_1$ then $\tilde g_t\in[0,3]$ 
for similar reasons. 
Thus the conditions of Lemma~\ref{doubleyous} are
satisfied.

For   polynomials or formal power series $P$ and $\hat{P}$, denote by
$P^+$ the  formal power series 
obtained by replacing all coefficients of $P$ by
their absolute values, and write $P\le \hat{P}$ if the
coefficient of any monomial in $P$ is no greater than the
corresponding coefficient in $\hat{P}$. We will use the obvious fact that
if
$P^+$ is absolutely convergent (for a particular assignment of the
indeterminates) then so is $P$.

With \eqn{wrecurse} in mind, and with the aim of obtaining the useful inequality~\eqn{goodone} below, define the power series
$ \gamma_t^*$ for each $t\in \smalls$ by
\bel{gammastar}
\gamma_t^*=
1+w_0^++\sum_{{u \in \smalls\atop h \in \functions}\atop (u,h)\ne
(t,\bz)} w(u,t,h)^+ \gamma^*_u
\prod_{i=1}^{k} \frac{1}{2-\gamma^*_{t_i}},
\ee
which by induction has a unique solution in formal power series with constant
terms all 1.  Then
$$
\frac{1}{2-\gamma^*_{t_i}}=\sum_{j\ge 0}
(\gamma^*_{t_i} - 1)^j
$$
and so by induction, all
coefficients of $\gamma^*_t$ are nonnegative for each $t\in \smalls$.
Thus 
$$
\frac{1}{2-\gamma^*_{t_i}}\ge \sum_{j\ge 0}
(1-\gamma^*_{t_i} )^j=\frac{1}{ \gamma^*_{t_i}}
$$
and, again by induction, comparing~\eqn{wrecurse}
with~\eqn{gammastar} gives
\bel{goodone}
\hatma_t^{\,+}\le
\gamma^*_t
\ee
    for each $t\in\smalls$.

Now consider summing the terms of $\gamma^*_t(n,p,\tilbfg)$ for $p$ and $f$ as in the
lemma, when $n$ is sufficiently large. Since all coefficients of $\gamma^*_t$ are nonnegative, we
are at liberty to sum the terms in any convenient order. It is immediate from the proof of
Lemma~\ref{doubleyous} that $w(u,t,h)^+=O(x)$ and $w_0^+=O(p)=O(x)$. It is now straightforward to
verify from~\eqn{gammastar}, by a sequence of successive approximations beginning with
$\gamma^*\approx 1$ for all $t$, that 
\bel{termbound}
\gamma^*_t(n,p,\tilbfg)=1+O(x).
\ee
The lemma now follows 
since from~\eqn{goodone}, and the fact that the constant terms in all $\hatma$'s and $\gamma^*$'s  are all 1, $(\hatma_t-1)^+\le\gamma^*_t-1$.
\qed
 If 
$p$ and $f$ satisfy the conditions of Lemma~\ref{converge}, we may treat
$\hatma_t(f)$ as a number, being the sum of the series, for $n$ sufficiently large.
Since we may ignore small values of $n$, and since $p$ is a function of $n$, this makes $\hatma_t(f)$ a real-valued function of $f$ and $n$, and henceforth in this section we treat it as such. 
\begin{prop}\lab{hatmaclose}
Uniformly for all $f\in\smallfunctions$ and $t\in\smalls \setminus\{t^*\}$,
$$
|\hatma_t(f)-\ga(f,t)|= O\left(\frac{x+\phi_t\lambda_\larges}{\lambda_t}\right)
$$
where $ \phi_{t}= n^{o(1)}$ for $t=t^*$ and  $\phi_{t}=1$ otherwise.
\end{prop}
\proof
An induction like the one proving Proposition~\ref{hammaclose} is used. 
The inductive hypothesis is
$$
|\hatma_t(f)-\hamma(f,t)|=
O_t\left(\frac{x+ \phi_t\lambda_\larges}{\lambda_t}\right),
$$
where $O_t$ denotes a bound depending only on $t$.
  The result then follows by Proposition~\ref{hammaclose}.

Suppose that
$f=\bz$. Then $h=\bz$ in \eqn{hammadef}
and the terms in \eqn{hatmadef} with $h\neq\bz$ are 0 because $\tilde g_{t_i}=0$ for all $i$ by~\eqn{hatmaf}. 
Hence, the products in~\eqn{hammadef} and~\eqn{hatmadef} are empty, and by simple (downwards) induction on $t$, $\hamma(f,t)=\hatma_t(f)$
for all $t\in\smalls$.

 It remains to prove the lemma when $f\neq\bz$, which we assume henceforth. 
 
Note that~\eqn{hatmadef} contains terms such that, for some values of $f$, the corresponding terms are excluded~\eqn{hammadef} because $f-h\notin\functions$.
For the inductive step, we  bound these  terms first. After this, we consider the
error caused by replacing
$\hamma(f-h,u)$ by
$\hatma_u(f)$ in (\ref{hammadef}), as well as $\hamma(f_i,t_i)$
by $\hatma_{t_i}(f)$,  and  $f_i(t_i)+1$ by
$f_i(t_i)$.

Since $\hatma_t(f)=1+O(x)$ by Lemma~\ref{converge},
  and  $\tilde{w}_0=O(p)$ and $\tilde{w}(u,t,h)=O(x)$  
from Lemma~\ref{doubleyous},
all  terms in the summation in~\eqn{wrecurse} are
$O(x)$.
If $f-h\not\in\functions$ in~\eqn{wrecurse}, so that
$f(t_{i^\prime})-h(t_{i^\prime})< 0$ for some $t_{i^\prime}$, then
$f(t_{i^\prime})=O(1)$ and so $g_{t_{i^\prime}}=O(1/\lambda_{t_{i^\prime}})$.
The contribution of such a term in (\ref{hatmadef}) is
$O\left(\lambda_u/\lambda_t\lambda_{t_{i^\prime}}\right)$, which in the case $t_{i^\prime}\prec u$ is $O(x/\mean_t)$.
On the  other hand, if $t_{i^\prime}=u$, we have the same situation as 
 in the second paragraph 
after~(\ref{inductsigma}), so the $\gamma$'s cancel, $c(u,t,h)=O(x)$, and again the term is $O(x/\mean_t)$.
 
For those $h$ satisfying $f-h\in\functions$,
first recall, as observed in
  the middle of the proof of Proposition~\ref{hammaclose}, 
the product
in~\eqn{hammadef}, which we will denote by $\Pi$, is
$O(1)$. Analogous to~\eqn{diff1} and~\eqn{diff2} in the proof of
Proposition~\ref{hammaclose}, we will show that, for the same values of $(u,h)$ as in that Proposition,
\eq\lab{diff3}
\left|\hamma(f-h,u)-\hatma_u(f)\right|\frac{\lambda_u}{\lambda_t}
c(u,t,h)\Pi
=
   \left\{ \begin{array}{ll}
O_u\left(\frac{x+\lambda_\larges}{\lambda_t}
\right)+O\left(\frac{x}{\mean_t}\right)&\mbox{
if }t\prec u\\
O_t\left(x\frac{x+ \phi_t\lambda_\larges}{\lambda_t}\right)
+O\left(\frac{x}{\mean_t}\right)&\mbox{ if }u=t,\end{array} \right.
\en
and 
\eq\lab{diff4}
\left|\hamma(f_i,t_i)-\hatma_{t_i}(f)\right|\frac{\lambda_u}{\lambda_t}
c(u,t,h)\Pi
=O_{t_i}\left(x\frac{x+\phi_{t_i}\lambda_\larges}{\lambda_t}\right)
+O\left(\frac{x}{\mean_t}\right),
\ \ \ t_i\prec u
\en
and, for the replacement of $f_i(t_i)+1$ by
$f_i(t_i)$ when evaluating $g_{t_i}$,
\bel{fichange}
\frac{\lambda_uc(u,t,h)}{\lambda_t\lambda_{t_i}}=
O\left(\frac{x}{\mean_t}\right).
\ee
The lemma follows from these claims, using~\eqn{deltas} along the lines of the proof of Proposition~\ref{hammaclose}, combined with the observation that, by the inductive hypothesis combinded with Lemma~\ref{converge}, we may assume that $ \hamma(f-h,u)=\Theta(1)$ uniformly whenever $h>\bz$, or  $h=\bz$ and $t\prec u$.

The treatment of the $O_t()$ terms in this proof is rather delicate and is explained in
detail in the proof of Proposition \ref{hammaclose}. In this case, there are extra terms
$O(x/\la_t)$ in (\ref{diff3}--\ref{fichange}), which we write separately to make the recursive argument clearer.
Note that the $O_u()$   and $O_t()$ terms contain the same implicit constants as in the
inductive hypothesis. 

It is convenient to treat~\eqn{fichange} first. If
$t_i\prec u$, then we are done by~\eqn{meanratio} applied with $t$ replaced by $t_i$, and the fact that
$c(u,t,h)=O(1)$. On the other hand, if $t_i=u$ then  $k=1$ and
$h=\delta_u$, and, 
as in the last part of the proof of Lemma~\ref{doubleyous}, 
$c(u,t,h)=O(x)$, as required.

 Now consider~\eqn{diff3}. Since either $f-h<f$ or $t\prec u$, the inductive hypothesis may be applied, with $\Pi$ referring to $f-h$ rather
than $f$, yielding
\begin{eqnarray}
|\hamma(f-h,u)-\hatma_u(f-h)|
\ut c(u,t,h) \Pi
&=&
O_u(1)
\frac{x+ \phi_u \lambda_\larges}{\lambda_u}
\ut c(u,t,h)\Pi\lab{xplusut}\\
&=&
O_u\left(
\frac{x+  \phi_u \lambda_\larges}{\lambda_t}c(u,t,h)\right). \nonumber
\end{eqnarray}
Recalling also from the proof of Lemma~\ref{doubleyous} that $c(t,t,h)=O(x)$ (and $c(u,t,h)=O(1)$
always),  and that $\phi_u=1$ when $t\prec u$, 
now shows that this expression
is bounded by $O_t\big( x(x+\phi_t\la_\larges)/\la_t\big)$ (respectively
$O_u\big((x+\la_\larges)/\la_t\big)$ ) as required for the cases $u=t$ and $t\prec u$ in the right hand
side of~\eqn{diff3}.
Next we bound
\bel{fminush}
|\hatma_u(f)-\hatma_u(f- h_0)| 
\ee
for any fixed $h_0$ with bounded entries. We can assume $h_0\ne \bz$. 
By Lemma~\ref{converge}, equation~\eqn{wrecurse} can be expanded in
increasing powers of the
$w$'s, which are
$O(x)$  under the substitution
$g_v=f(v)/\lambda_v$ by Lemma~\ref{doubleyous}. By~\eqn{xbound}, we may ignore
terms whose total degree in $w$'s is larger than some fixed value. Into the truncated
expression, substitute
$f(t_i)/\la_{t_i}$ and  $(f(t_i)-h(t_i))/\la_{t_i}$ for $g_{t_i}$ in the
definition of
$w(u,t,h)$ at  (\ref{wdefs})  and subtract
the two resulting expressions term by
term.     Since the entries of $h_0$ are bounded, the  dominating terms are exactly of the type estimated
in~\eqn{fichange}, and hence are bounded by $O\left(x/\mean_t\right)$. Equation~\eqn{diff3} now follows  (with room to spare) 
in view of the fact that,  by Lemma~\ref{doubleyous},
$$
 \ut c(u,t,h)\Pi=O(x).
$$

   The proof of (\ref{diff4}) involves  firstly consideration of
$\left|\hamma(f_i,t_i)-\hatma_{t_i}(f_i)\right|$ (multiplied by the other
factors).
This yields an expression as in
the right hand side of~\eqn{xplusut}, but with
$O_{u}$ replaced by $O_{t_i}$, $\lambda_\larges/\lambda_u$ replaced with
$\lambda_\larges/\lambda_{t_i}$ and $f-h$ becoming $f_i$. The error term is bounded similarly to the bound
 {\eqn{tibound}} for the analogous term in the proof of Proposition~\ref{hammaclose},
 and also using  $\lambda_u/\lambda_{t_i}=O(x)$ (as $t_i\prec u$),
giving the first error term in~\eqn{diff4}.
  Then, $\left|\hatma_{t_i}(f_i)-\hatma_{t_i}(f)\right|$ is bounded by the
  expression in~\eqn{fichange}, by the same argument as for~\eqn{fminush}.
\qed

{From} Lemmas~\ref{doubleyous} and~\ref{converge}, we may use~\eqn{wrecurse} to expand all the
functions $\hatma_t$ ($t\in \smalls$) recursively in power series in $n$, $p$
and the variables $g_i$. Iterating $r$ times determines $\hatma_t$ to arbitrarily small error
$O(x^r)$ when the appropriate values are assigned to $p$ and the $g_i$. However, instead of pursuing arbitrary accuracy in this paper, we desire a final formula which is shown to exhibit a uniformity over all relevant $\kappa$, and for this we need the following.   We use $\bfg^{\bf i}$ to denote $g_1^{i_1}g_2^{i_2}\cdots g_s^{i_s}$;  if ${\bf i}={\bf 0}$, this is the multiplicative identity   of the ring $\real[[n,p]][[\bfg]]$ of formal power series over $\bfg$ whose coefficients are in $\real[[n,p]]$.

\begin{cor} \lab{expanded}
There are power series $\xi_t$, $t\in\types$, in $n$, $p$, $\bfg$, independent of  $\kappa$, and, for all $\eps>0$, truncations $\xi_{t,\eps}$ of the
series   $\xi_t$, to a finite number of terms, such that for all $t\in\smalls$
\begin{description}

\item{(a)}  For   ${\bf i}\ne  {\bf 0}$ , we have $[\bfg^{\bf i}]\xi_t= O(x)$,
 for $p$ satisfying \eqn{pkappa} with  $\kappa\ge\chi+\eps$,
 as $n \to \infty$;
\item{(b)}  for each {\bf i}, the   coefficient $[\bfg^{\bf i}]\xi_t$ is a multiple of $\prod_{u\in\smalls} p^{\mu(u)i_u}$;
\item{(c)} With $p$ satisfying~\eqn{pkappa},
and  $\xi_{t,\eps}(f)$ defined from  $\xi_{t,\eps}$ analogously to $\hatma_t(f)$ in~\eqn{hatmaf},
there exists $\epshat>0$ such that uniformly  for all $f\in \smallfunctions$, and   all $\kappa\ge\chi+\eps$, 
\bel{excor}
\xi_{t,\eps}(f)=\hatma_t(f)+O\left(\frac{ n^{-\epshat}}{\lambda_t}\right).
\ee 

\end{description}
\end{cor}

\proof 
Instead of (c) we show the obviously stronger
\bel{excormain}
\xi_{t,\eps}(f)=\hatma_t(f)+O\left(\frac{x+n^{o(1)}\lambda_\larges}{\lambda_t}\right). 
\ee 
 We start by essentially focusing on this, but with one eye fixed on (a). Define the function $F_t=F_t(n,p,\bfg, \hatma_{ 1}, \ldots , \hatma_{ s})$ by
\eq\lab{Ftdef}
F_t(n,p,\bfg, \hatma_{ 1}, \ldots , \hatma_{ s})=
\frac{1}{c(t,t,\bz)}\left(
1-\sum_{{u \in \smalls\atop h \in \functions}\atop (u,h)\ne (t,\bz)}
\frac{\lambda_u}{\lambda_t}c(u,t,h)\hatma_u
\prod_{i=1}^{k} \frac{g_{t_i}}{\hatma_{t_i}}
\right)-1.
\en
We obtain successive power series approximations $F_t^{(j)}$ and $\hatma^{(j)}_t$ for all the $F_t$ and $\hatma_t$ ($j=0,1,\ldots$).
Initially,    set  $F_t^{(0)}=0$ and
$\hatma^{(0)}_t = 1$
for all $t$. For $j\ge 0$, substituting $\hatma_t^{(j)}$ for $\hatma_t$
in (\ref{Ftdef}) simultaneously
for all $t\in\smalls$ defines $F_t^{(j+1)}$ as a power series (recalling the observations made before~\eqn{hatmadef} that $\la_u/\la_t$ is a polynomial in $n$ and $p$, and so on). Next, define  $\hatma^{(j+1)}_t= 1+  F_t^{(j+1)}$ to complete the iterative definition.
Define  $\hatma^{(i)}_t(f)$ from
$\hatma^{(i)}_t$ analogously to
$\hatma_t(f)$ in~\eqn{hatmaf}, and similarly $F_t^{(i)}(f)$.  By
Lemma~\ref{converge},
$\hatma^{(0)}_t(f)=
\hatma_t(f)(1+O(x))$  for all relevant $f$ and $p$. Thus
$$
F_t^{(1)}(f)= F_t(n,p,\tilbfg, \hatma_{ 1}, \ldots , \hatma_{ s})(1+O(x)).
$$
By Lemma~\ref{converge}, this is $O(x)$, and so by~\eqn{hatmadef}, $\hatma^{(1)}_t(f)=\hatma_t(f) +O(x^2)$.  Repeating the same argument $r$
times shows that
\bel{hatmar}
\hatma^{(r)}_t(f)=\hatma_t(f) +O(x^{r+1}).
\ee

As with Lemma~\ref{converge}, the argument to this point is for fixed $\kappa>\chi$.
The definition of $\smalls$ by (\ref{smallsdef}), and hence the
formula~\eqn{hatmadef}, depends on $\kappa$. However, for all $\kappa\ge\chi+\eps$,
$\smalls$ is a subset of
$\hat \smalls=\{t\in\types:\nu(t)/\mu(t)\ge \chi+\eps\}$,  which is the value of
$\smalls$ when $\kappa= \kappa_0 = \chi+\eps$. So define
$r_t$ to be such that   $x^{r_t}=O(1/\la_t)$ when $\kappa =\kappa_0$. Then  set
$\xi_{t,\eps}$ equal to the truncation of $ \hatma^{(r_t)}_t$ to those terms whose
value, with $\bfg$ set equal to 1,  is not $o(x/\la_t)$ (when $\kappa=\kappa_0$).
By~\eqn{hatmar},~\eqn{excormain} holds for $\kappa=\kappa_0$. 

  Also note for later use
that, in view of~\eqn{hatmar}, using
$ \hatma^{(r)}_t$ for any $r>r_t$ would define the same $\xi_{t,\eps}$.
  From~\eqn{goodone} and~\eqn{termbound}, the  coefficients of any non-constant monomial   $\bfg^{\bf i}$ in $\xi_{t,\eps}$, as it arises recursively from~\eqn{Ftdef}, are $O(x)$, which proves part (a) with $\xi_{t}$ interpreted as $\xi_{t,\eps}$.

We next claim that~\eqn{excormain} is also valid when $\kappa>\kappa_0$. In this case,
the recursive definition of $\hatma^{(r)}_t$ is the same as for $\kappa_0$ except
that the definition of $\smalls$ is different. Any terms in the summation
in~\eqn{hatmadef}  corresponding to types
$t$ that are in $\smalls$ for $\kappa_0$, and not in $\smalls$ for $\kappa$, are now
missing.  These terms are of the form $\la_u/\la_t$ times a finite product of $g_i$,
for some
$u\notin\smalls$. Since all $g_i$ are substituted with values $n^{o(1)}$, the claim
holds.

The remaining portion of the claim in part (c) of the corollary relates to uniformity. This follows from the above observations once we   show that these functions $\xi_{t,\eps}$ are all common
truncations of the power series $\xi_{t}$. Now of course (a) is justified in its original form, for $\xi_{t}$.

If $\eps'<\eps$ is considered,
then new types enter $\smalls$, but the terms in $\xi_{t,\eps'}$ due to these are
of smaller order (as with consideration of $\kappa>\kappa_0$ above) and cannot be
included in $\xi_{t,\eps}$. Also, the appropriate value of $r_t$ may be larger for
$ \eps' $ than for $\eps$, but as noted above, truncating with the larger value of
$r$ gives the same function $\xi_{t,\eps}$, so the extra terms generated cannot
include any of the same monomials as appearing in $\xi_{t,\eps}$.
The power series $\xi_t$ is now well-defined to be the termwise limit of  
$\xi_{t,\epsilon}$ as $\epsilon\to 0$.

Finally, to verify part (b), note that in the recursive use of~\eqn{Ftdef}, every new product
$\prod_{i=1}^{k}  g_{t_i}$ that is  introduced is accompanied by the factor $\frac{\lambda_u}{\lambda_t}c(u,t,h)$. By its definition~\eqn{cdef}, each term of $c(u,t,h)$ is associated with a cluster of $J$ of type $u$, a cluster $Q$ of type $t$, and pairwise edge-disjoint clusters $J_1,\ldots, J_k$ of types $t_1,\ldots, t_k$, with  $c(u,t,h)$ divisible by $p^a $ where $a=\left|Q
\cap \left(\bigcup J_i\right)\right|$. Since $\la_u/\la_t$ is divisible by $p^b$ where $b= \mu(u)-\mu(t) = \mu(u) - |Q|$, the term itself must be divisible by $p^{\sum|J_i|}$, as required for part (b).
Of course, the expansions of $1/c(t,t,\bz)$ and $1/\hatma_{t_i}$ do not affect this as their terms have nonnegative exponents.
\qed


\section{Graphs with forbidden subgraphs in $\Gnp$}\lab{s:main}

In this section we prove our main result for subgraphs of the random  graph $\Gnp$.
Let $G_0$ be a strictly balanced graph and recall that $\chi$ is defined by
(\ref{chidef}). Let $X$ be the number of copies of $G_0$ in $\Gnp$.  
  \bigskip

\noindent
{\bf Proof of the $\Gnp$ case of Theorem~\ref{t:main}}

The proof works roughly as follows. We estimate the ratios of `adjacent' probabilities $\pr(\class_f)$ by estimating $\ga(f,t)$ defined in~\eqn{gammadef}. This is approximated by $\hatma_t(f)$, as shown in Proposition~\ref{hatmaclose}, which in turn is approximated by $\xi_{t,\eps}$ as found in Corollary~\ref{expanded}.
 Fix $\epsilon>0$. We   assume at first that $p=n^{-\kappa+o(1)}$ for fixed $\kappa\ge\chi+\eps$,
in accordance with~\eqn{pkappa}, so that~\eqn{xbound}, Proposition~\ref{hatmaclose} and Corollary~\ref{expanded} can be applied.  The theorem will then be shown in full
generality, with assistance from Lemma~\ref{panything}. In this section, we work only with the $G_0$-clustering.  As a consequence of this, the parts of the theorems in the previous section relating to $t^*$ are not needed.   The set $\smalls$ is defined, as before, to contain just
those types $t$ in this clustering for which $\nu(t)/\mu(t)\geq\kappa$. Recall by the discussion after (\ref{pkappa})
that $\smalls$ is finite.

The expected number of sets of $j$ disjoint clusters of type $t\in\smalls$ is,
recalling~\eqn{tsize} and~\eqn{lambdadef}, at most
$$
 {|t|\choose j}p^{\m(t)j}
 \leq
\left(\frac{e|t|p^{\m(t)}  }{j }\right)^j = \left(\frac{e \lambda_t}{j}\right)^j.
$$
Taking $j=\lfloor m_t\rfloor+1$    for each $t\in\smalls$ shows by~\eqn{mtdef} (using $e<3$) that
$\sum_{f\not\in\smallfunctions}\BP(\class_f) =o(1)$. (This reveals the relevance of the constant 3
in the definition of $m_t$.) Furthermore, every large cluster contains an unavoidable cluster, of
which there are a finite number. Applying~\eqn{maxlargebound} to all such clusters, we see that
$\sum_{f\in\functions}\BP(\class_f) \sim 1$. Hence 
\eq\lab{large} 
\BP(X=0)^{-1}=\frac{1}{\BP({\cal C}_\bz)}\sim
\sum_{f\in\smallfunctions} \frac{\BP({\cal C}_f)}{\BP({\cal C}_\bz)}. 
\en

By  renaming the cluster types in $\smalls$ if necessary, 
extend the poset on $\smalls$ to a unique linear ordering
on $\smalls=[s]:=\{1,2,\ldots,s\}$ denoted by $<$, in  decreasing order of 
 $\nu(t)-\kappa \mu(t)$, breaking ties in a canonical way independent of the choice of $\kappa$ (i.e.\ depending only on the graph structure of the types).  This is possible in view of~\eqn{tsize},~\eqn{lambdadef}, and~\eqn{meanratio}. Although the values of $p$ can ``wobble'' around $p^{-\kappa}$, so that $\la(t+1)$ and $\la(t)$ are not  always in the same order when a tie occurred, we do have
\bel{labound}
\la_{t+1}<n^{o(1)}\la_t \quad\mbox{ for all $t<s$}. \ee
 (That observation is in fact  the main   motivation behind the restriction of $p$ in (\ref{pkappa}).)

 Fix $ (j_1,j_2,\ldots,j_{s})$ with $j_u\in [0,m_u]$ for all $u\in \smalls$ and define $f$ so that $f(t)=j_t$ for each $t\in \smalls$. Then for each $t$ and $j$ define the function
$f_{t,j}$ on $\smalls$ by $f_{t,j}(t^\prime)=j_{t^\prime}$ for $t^\prime<t$; $f_{t,j}(t)=j$;
$f_{t,j}(t^\prime)=0$ for $t^\prime>t$. Then $ f_{s,j_s}=f$. 

 By
Proposition~\ref{hatmaclose}, we have $\hatma_t(f)=\gamma(f,t) + O((x+\mean_\larges)/\lambda_t)$
uniformly for all $f\in\smallfunctions$ and $t\in\smalls$. Moreover, by
Proposition~\ref{thetabound2}, $\gamma(f,t)\sim 1$ uniformly, so that
$\hatma_t(f)=\gamma(f,t)\big(1+O((x+\mean_\larges)/\lambda_t)\big)$. Note that   $\big(1+O((x+\mean_\larges)/\lambda_t)\big)  = \big(1+O(  n^{-\epshat}/\lambda_t  )\big) $   by~\eqn{xbound} and~\eqn{maxlargebound}. Using these estimates, then Corollary~\ref{expanded}, and finally the fact that 
$\big(1+O(  n^{-\epshat}/\lambda_t  )\big)^{m_t}=1+o(1)$ by the definition of $m_t$
in~\eqn{mtdef}, we have  
\bea
\frac{  \pr(\class_f)}{\pr(\class_{\bf 0})} &=&
\prod_{t=1}^s\prod_{j=0}^{j_t-1}\rho(f_{t,j},\delta_t)\nonumber \\
&=&
 \prod_{t=1}^s\frac{\lambda_t^{j_t}}{j_t!} \prod_{j=0}^{j_t-1}\gamma\left(f_{t,j},t\right) \nonumber\\
&=&
\ \prod_{t=1}^s\frac{\lambda_t^{j_t}}{j_t!} \prod_{j=0}^{j_t-1}\hatma_{t}\left(f_{t,j}\right)\big(1+O((x+\mean_\larges)/\lambda_t)\big) \lab{hatmaest}\\
&=&
 \prod_{t=1}^s\frac{\lambda_t^{j_t}}{j_t!} \prod_{j=0}^{j_t-1}\xi_{t,\eps}\left(f_{t,j}\right)\big(1+O(  n^{-\epshat}/\lambda_t  )\big)
\nonumber\\
&=&
\big(1+o(1)\big) \prod_{t=1}^s\frac{\lambda_t^{j_t}}{j_t!} \prod_{j=0}^{j_t-1}\xi_{t,\eps}\left(f_{t,j}\right).\lab{pointprob}
\eea

Our basic method is to sum the above expression over all $f$ for which $\pr(\class_f)$ is significant, thereby obtaining an estimate for the reciprocal of $\pr(\class_{\bf 0})$. To facilitate analysis of the summation, we employ various partial sums defined as follows. For $t\in \smalls\cup \{0\}$, define the functions $S_t$ by
$S_s(j_1,j_2,\ldots,j_s)=1$, and recursively  for  $t$ decreasing from $s-1$ to 0, by
\eq\lab{Sdef}
S_t(j_1,j_2,\ldots,j_t) = \sum_{j=0}^{\lfloor m_{t+1}\rfloor } S_{t+1}(j_1,j_2,\ldots,j_t,j)
\left( \prod_{i=0}^{j-1} \hatma_{t+1}\left(f_{t+1,i}\right) \right) \frac{\lambda_{t+1}^j}{j!}.
\en 
We next show (see~\eqn{Srecursion}) that this quantity approximates the reciprocal of the conditional probability of having no small
clusters of type $u>t$, given $j_u$ clusters of type $u$ for all $u\leq t$. Recalling the error bounds involved   in~\eqn{hatmaest}, and then the bound on $m_t$ used in deriving~\eqn{pointprob},  we have, uniformly,
\begin{eqnarray*}
S_t(j_1,j_2,\ldots,j_t)
&\sim& \sum_{j=0}^{\lfloor m_{t+1}\rfloor} S_{t+1}(j_1,j_2,\ldots,j_t,j) \left( \prod_{i=0}^{j-1}
\gamma\left(f_{t+1,i},t+1\right) \right)
\frac{\lambda_{t+1}^j}{j!}\\
&=& \sum_{j_{t+1}=0}^{\lfloor m_{t+1} \rfloor } S_{t+1}(j_1,j_2,\ldots,j_t,j_{t+1})
\frac{\pr\left(\class_{f_{t+1,j_{t+1}}}\right)} {\pr\left(\class_{f_{t+1,0}}\right)}.
\end{eqnarray*}
An inductive argument immediately shows that for all $t<s$ and $j_i\in[0,m_t]$, $i\in[1,t]$,
\eq\lab{Srecursion} S_t(j_1,j_2,\ldots,j_t) \sim \sum_{j_{t+1}=0}^{\lfloor m_{t+1}\rfloor } \cdots
\sum_{j_s=0}^{\lfloor  m_s \rfloor} \frac { \BP( \class_{f_{t+1,0}+ j_{t+1}\delta_{t+1}
+\cdots+j_s\delta_s })} {\BP(\class_{f_{t+1,0}})}. \en Therefore, by~\eqn{large}, noting that
$S_0$ has no arguments, \eq\lab{S0sim} S_0\sim \BP(X=0)^{-1}. \en Thus, we  have reduced the
problem to that of estimating $S_0$.

For use in the following, we define \eq\lab{hatgdef} \zeta(j ) =\left( j_1 /\lambda_1,  \ldots,
j_{t-1} /\lambda_{t-1}, j /\lambda_t\right) \en (with the dependence on $j_1,\ldots ,j_{t-1}$
suppressed for compactness of notation), and we say that $\zeta(j)$ is {\em appropriate} if
$j_i\in[0,m_i]$ for all $i\in[1,t-1]$ and $j\in [0,m_t]$.

It is useful to define $\calP^+(\bfgt)$ to be the ring of polynomials in $\bfgt = (g_1,\ldots, g_t)$
whose  coefficients are polynomials in $n$, $p$ and $n^{-1}$ (as formal indeterminates), and $\calP(\bfgt)$ to be the subring of $\calP^+(\bfgt)$ consisting of those polynomials whose coefficient of
$g_1^{i_1}\cdots g_t^{i_t}$ is divisible by
 $p^{\sum_j \mu(j)i_j}$.   We regard these coefficients simply the union of the set of  polynomials in $n$ and $p$   with the set of polynomials in $n^{-1}$ and $p$, which together form a ring.

We will use the definition of $S_t$, together with Lemma~\ref{converge} and
Corollary~\ref{expanded} and an induction argument to prove that
\eq\lab{Sform}
S_0=\exp \big( P_{0,\kappa} +o(1) \big), \qquad
S_t(j_1,\ldots,j_t)=\exp \big( P_{t,\kappa} \big(\zeta( j_t)\big)+o(1) \big) \quad
(1\le t \le s) 
\en
  for all $j_1,\ldots,j_t$ such that $\zeta(j_t)$ is appropriate, for some polynomials  $P_{t,\kappa}$ such that
 \begin{description}
\item{(i)}
  $ P_{t,\kappa}\in \calP(\bfgt)$   (and so in particular for $t=0$, $P_{0,\kappa}$ is a polynomial in $n$, $p$ and $n^{-1}$);
\item{(ii)}
the constant coefficient of $P_{t,\kappa}$ (i.e.\ $P_{t,\kappa}(0,0,\ldots,0)$)  is
 equal to $(1+O(xn^{o(1)}))\sum_{t'=t+1}^s\la_{t'}$
and the other coefficients are $O(xn^{o(1)}P_{t,\kappa}(0,0,\ldots,0) )$, where the implicit
bounds  in $O(\cdot)$  are  independent of $\eps$. 
Note that by~\eqn{labound}, it follows that the constant coefficient is
$O(n^{o(1)}\la_{t})$;
\item{(iii)}  the convergence expressed by $o(1)$ in~\eqn{Sform} is uniform over all appropriate $j_1,j_2,\ldots,j_{t}$.
\end{description}
The induction begins with $t=s$ and then proceeds through decreasing values of $t$. It finishes with the case $t=0$ of (\ref{Sform}), which is used to show that the polynomials $P_{t,\kappa}(\bfgt)$ are of
such a form that the theorem follows using (\ref{S0sim}).

The initial step of the induction argument, $t=s$, is trivial, since $S_s$ is identically equal to
1 and we may set $P_{s,\kappa}=0$. So now suppose   that (\ref{Sform}) holds for some particular value of $t$. We must prove that
it also holds when $t-1$ is substituted for $t$. Define $T_j$ by
\bel{Tdef} 
T_{j} =\exp \big(
P_{t,\kappa} \big(\zeta(j ) \big)\big) \left(\prod_{i=0}^{j-1}
\xi_{t,\eps}\left(f_{t,i}\right)\right) \frac{\lambda_t^j}{j!}. 
\ee 
 We now use~\eqn{Sdef},~\eqn{Sform} and   Corollary~\ref{expanded}  to replace $\hatma$ in~\eqn{Sdef} by $\xi$, the fact that  $\zeta(j)$ is appropriate and  $m_t=O(\la_t \log n)$, together with~\eqn{maxlargebound}, to obtain
\bel{Stm1d} 
S_{t-1}(j_1,j_2,\ldots,j_{t-1}) \sim \sum_{j=0}^{\lfloor  m_{t} \rfloor} T_j. 
\ee

First assume that $t\in\smalls_0$. Note that (with square brackets for extraction of coefficients)
$$
\frac{\exp  P_{t,\kappa} \big(\zeta(j)\big)}
{\exp   P_{t,\kappa} \big(\zeta(0)\big)} =
\exp \sum_{{\bf i}} \big([ \bfg^{\bf i}]P_{t,\kappa}\big)
\left(\prod_{\ell=1}^{t-1} (j_\ell/\la_\ell)^{i_\ell}\right)
\big((j/\la_t)^{i_t}-0^{i_t}\big) 
$$
where $0^0=1$ as usual, and the summation is over the set of $\bf i$ for which the coefficient is nonzero. The number of such $\bf i$ is bounded, given $P_{t,\kappa}$. 
The only terms contributing have $i_t>0$, and in particular the constant term does not contribute.  Let $\vmax$ be the total degree of $P_{t,\kappa}(\bfgt)$. Each factor $j_\ell/\la_\ell$ is by~\eqn{s0s1} at most  $\log n =  n^{o(1)}$, and the same goes for $j/\la_t$.
By the inductive hypothesis (ii) and~\eqn{mtdef}, we now obtain
$$
\frac{\exp  P_{t,\kappa} \big(\zeta(j)\big)}
{\exp   P_{t,\kappa} \big(\zeta(0)\big)} =
\exp\big(xn^{o(1)}(n^{o(1)})^\vmax\big)
=1+O(xn^{o(1)}).
$$
By  Lemma~\ref{converge}, $\hatma(t)=1+O(x)$, and so Corollary~\ref{expanded}  gives that   each factor $\xi_{t,\eps}\left(f_{t,i}\right)$  in~\eqn{Tdef}  is 
$1+O\big(x+n^{-\epshat+o(1)}\big)$
Hence for $j \le m_t =n^{o(1)}$,  the product of $j$
factors in~\eqn{Tdef}  is
$$
{\left(1+O\big(x+n^{o(1)-\eps'}\big)\right)}^{m_t}\sim 1
$$
using~\eqn{xbound}. Thus
$$
T_j \sim  \exp   P_{t,\kappa} \big(\zeta(0)\big) \lambda_t^j/j! \sim S_t(j_1,j_2,\ldots,j_{t-1},0)
\lambda_t^j/j! 
$$
by the inductive hypothesis~\eqn{Sform}.
Since in this case $m_t=\la_t \log n$, it follows that
$$
S_{t-1}(j_1,j_2,\ldots,j_{t-1}) = S_t(j_1,j_2,\ldots,j_{t-1},0) \exp\left(\lambda_t+o(1)\right).
$$
Here we used the uniformity of the convergence in the estimates, including that asserted in part
(iii) of the induction hypothesis. To  establish  the inductive hypothesis in this case, we thus
set $P_{t-1,\kappa}$   equal to $P_{t,\kappa}+\la_t$, which is a polynomial in $n$ and $p$ (thus a constant in
$\calP(\bfg)$) by~\eqn{tsize} and~\eqn{lambdadef}. This clearly gives the inductive
hypotheses (i) and (ii), whilst the uniformity in (iii) implies that (iii) holds with $t$ replaced by $t-1$.

We next suppose   that $ t\in \smalls_1$, so that in particular
  $\lambda_t\to\infty$ by~\eqn{s0s1}.
We need to estimate the ratio of consecutive terms $T_j$ quite accurately. 
We have 
\bel{ratio} 
\frac{T_j}{T_{j-1}} = \exp  \Big(P_{t,\kappa}
\big(\zeta(j)\big)-P_{t,\kappa} \big(\zeta(j-1)\big)\Big) \xi_{t,\eps}
(f_{t,{j-1}})\frac{\la_t}{j}.
\ee

Let  $R_v=[g_t^v]P_{t,\kappa}(\bfgt)$, so that $R_v\in\calP(\bfg_{t-1})$. Put
$$
\hat\zeta = (j_1/\la_1,\ldots, j_{t-1}/\la_{t-1})
$$
and  
\bel{etadef} 
\eta = n^{-\eps/2} /\la_t.
\ee
Then
\bean
 P_{t,\kappa} \big(\zeta(j)\big)-P_{t,\kappa} \big(\zeta(j-1)\big)&=&
 \sum_{v=1}^{v_{\rm max}} R_v(\hat\zeta )\Bigg(\bigg(\frac{j}{\la_t}\bigg)^v- \bigg(\frac{j-1}{\la_t}\bigg)^v\Bigg)\\
& =& \sum_{v=1}^{v_{\rm max}} R_v(\hat\zeta )  \left(\frac{vj^{v-1}+O(j^{v-2})}{\la_t^{v}}\right) \\
& =& \sum_{v=1}^{v_{\rm max}}\frac{vR_v(\hat\zeta )}{\la_t} \cdot \frac{ j^{v-1}}{\la_t^{v-1}} + O (R_v(\hat\zeta )  /\la_t^2)\\
& =&
  O(\eta )
 +  \sum_{v=1}^{v_{\rm max}}\frac{vR_v(\hat\zeta )}{\la_t} \cdot \frac{ j^{v-1}}{\la_t^{v-1}}
\eean since $j=O(\la_t)$ by~\eqn{mtdef} (and $j_i=O(\la_i)$ for $i<t$), and using the inductive
hypothesis (ii), which implies that   the coefficients of $R_v$ for $v\ge 1$  are all
$O(n^{o(1)}x\la_t)=O(\eta\la_t^2)$ by~\eqn{xbound}. For the same reason, the  terms in this
summation are  all $O(n^{o(1)}x)$.

We  call a polynomial $\tilde{P}$ {\em acceptable} if
$\tilde{P}=1+P$   for some polynomial $P\in\calP(\bfg)$ whose   coefficients are all
$O(n^{o(1)}x)$ for the range of $p$ under consideration, i.e.\ satisfying~\eqn{pkappa}. Note that
$ n^{o(1)}x=o(n^{-\eps/2})=o(1)$ by~\eqn{xbound}. 
A polynomial $\tilde{P}$ is {\em $t$-acceptable} if  $\tilde{P}=1+P$   for some polynomial $P\in\calP^+(\bfgt)$  whose coefficient of
$g_1^{i_1}\cdots g_{t}^{i_{t}}$ is divisible by
 $p^{\sum_{j<t} \mu(j)i_j}$, and whose   coefficients are all
$O(n^{o(1)}x)$ for $p$ satisfying~\eqn{pkappa}.  That is,
 $\tilde{P}$ satisfies the definition of an acceptable polynomial in $\calP(\bfgt)$ except that the powers of $p$ in the terms in $P$ are only required to pay their respect to the variables $g_1,\ldots, g_{t-1}$.  

By~\eqn{lambdadef}, $\la_t^{-1}$ can be
expanded as $p^{-\mu(t)}$ times a power series in $n^{-1}$. So by the inductive assumption that
$P_{t,\kappa}\in \calP(\bfgt)$, it follows that there exists $\tilde{R}_v \in \calP(\bfg_{t-1})$ such
that $R_v(\hat{\zeta})/\la_t=\tilde{R}_v(\hat{\zeta})+O(\eta)$. To verify this, we note that $g_t$ does not appear
in $R_v$ and hence the lower bound on the exponent of $p$ required for $P_{t,\kappa}$'s membership in $\calP(\bfgt)$
is enough to compensate for $p^{-\mu(t)}$; the power series in $n^{-1}$ can be truncated at an
appropriate point to obtain a polynomial in $n^{-1}$, producing the error term $O(\eta)$.

We conclude that 
$$
 \exp\Big( P_{t,\kappa} \big(\zeta(j)\big)-P_{t,\kappa} \big(\zeta(j-1)\big)\Big) =   A_{t,\kappa}^{(1)} (\zeta(j))(1+O(\eta))
 $$
for a $t$-acceptable polynomial $A_{t,\kappa}^{(1)}$ (with constant term precisely 1 in this case).

By Lemma~\ref{converge} and Corollary~\ref{expanded}(b), we see that  $\xi_{t,\eps} $ is
acceptable and consequently $t$-acceptable. Consequently,  $\xi_{t,\eps}(f_{t,{j-1}})$ that occurs in~\eqn{ratio} is equal to $\xi_{t,\eps}(f_{t,{j}})(1+O(\eta))$. Moreover, the product of two $t$-acceptable polynomials is $t$-acceptable. Thus~\eqn{ratio}
gives
 \bel{ratio2}
\frac{T_j}{T_{j-1}}=\frac{ A_{t,\kappa}^{(2)} (\zeta(j))\la_t(1+O(\eta))}{j} \ee for the  $t$-acceptable
polynomial 
\bel{A2def}
A_{t,\kappa}^{(2)} (\bfgt):= A_{t,\kappa}^{(1)}\cdot \tilde\xi_{t,\eps},
\ee
 where $\tilde\xi_{t,\eps}$ is obtained from $ \xi_{t,\eps}$ by setting $g_{t+1}=\cdots =g_s=0$.

 To identify (approximately) the maximum term of the summation in~\eqn{Stm1d}, we  note that since $A_{t,\kappa}^{(2)}$ is $t$-acceptable,  $A_{t,\kappa}^{(2)} (\zeta(j))\sim 1$ and so~\eqn{ratio2} shows that we are interested in $j\sim \la_t$. Furthermore, again using $t$-acceptability, the derivative of $A_{t,\kappa}^{(2)}(g_1,\ldots, g_{t-1},y)$  with respect to $y$   is $o(n^{-\eps/2})$ when $\bfg_{t-1}$ is set equal to $\hat \zeta$. So, at least for large $n$, this function has a fixed point $y$ that is $1+o(1)$. In other words,  there must exist $j^*\sim \la_t$ satisfying
\bel{jstardef} 
j^* = \la_t   A_{t,\kappa}^{(2)}(\zeta(j^*)). 
\ee 
Since $A_{t,\kappa}^{(2)}$ is
$t$-acceptable, we can use repeated substitutions in
$$
Q_\ell =   A_{t,\kappa}^{(2)}(g_1,\ldots , g_{t-1},Q_{\ell-1})
$$
beginning with $Q_0=1$ to obtain a polynomial $Q_\ell\in  \calP(\bfg_{t-1})$ such that $Q_\ell(\hat \zeta)$ is an approximation to $j^*/\la_t$.   Clearly, replacing
the variable $g_t$  
of a  $t$-acceptable polynomial by another $t$-acceptable polynomial produces yet
another $t$-acceptable polynomial. So each $Q_\ell$ is $t$-acceptable.
For each iteration, the error in the approximation is multiplied by $o(n^{-\eps/2})$. Hence, for $\ell$ sufficiently large, $Q_\ell$ is an  acceptable polynomial 
$ A_{t,\kappa}^{(3)}\in \calP (\bfg_{t-1})$ satisfying
\bel{jstarasymp} 
j^*=\la_t A_{t,\kappa}^{(3)}(\hat\zeta )+o(1) 
\ee 
uniformly for all $\hat\zeta $ under consideration.
 
 For the product in~\eqn{Tdef} we will use the following.
Recall that $\xi_{t,\eps}$ is a polynomial, whereas $\xi_{t,\eps}\left(f_{t,i}\right)$ is a number given $n$ and $p$ (and in the present context $n$ determines $p$). Since $\xi_{t,\eps}$ is acceptable, we may expand its logarithm and hence obtain
\bel{logxi}
 \log\xi_{t,\eps}\left(f_{t,i}\right)= \sum_{v=0}^{\vmaxx}  R_v^{(1)}(\hat\zeta ){\left(\frac{i}{\lambda_t}\right)}^v
+o(\la_t^{-1})
 \ee
for some $\vmaxx$,  with $ R_v^{(1)}\in \calP(\bfg_{t-1})$ having   all coefficients  
$O(n^{o(1)}x)$ for all $v\le \vmax$. (That is, $1+  R_v^{(1)}$ is acceptable.) Then
\begin{eqnarray}\lab{xiexp}
\sum_{i=0}^{j-1}\log\xi_{t,\eps}\left(f_{t,i}\right) &=& \sum_{i=0}^{j-1}\sum_{v=0}^{\vmaxx}
R_v^{(1)}(\hat\zeta )  {\left(\frac{i}{\lambda_t}\right)}^v +
o\left(j/\lambda_t\right)\nonumber\\
&=& \sum_{v=0}^{\vmaxx} \frac{R_v^{(1)}(\hat\zeta )}{v+1}\cdot \frac{j^{v+1}}{\lambda_t^v} +
\sum_{v=0}^{\vmaxx}   \frac{O(j^v)R_v^{(1)}(\hat\zeta )}{\lambda_t^v}+o\left(
j/\lambda_t\right)\nonumber\\
&=&  o(1) + \lambda_t \sum_{v=0}^{\vmaxx} \frac{R_v^{(1)}(\hat\zeta )}{v+1}
{\left(\frac{j}{\lambda_t}\right)}^{v+1}.
\end{eqnarray}

We wish to approximate the terms in~\eqn{Stm1d} by expanding the formula for  $T_j$ given in~\eqn{Tdef}   about $j= j^* $, beginning with~\eqn{ratio2} written as
 \bel{Stexp}
 \log (T_j/T_{j-1}) = q(j)  +O(\eta)
 \ee
 where
 \bel{qdef}
 q(j)=\log A_{t,\kappa}^{(2)} (\zeta(j))+ \log \lambda_t    -\log j.
 \ee
Note that this equation also defines $q(y)$ for an arbitrary non-integer real $y$, so we can
consider its derivative $q'(y)$. Since $A_{t,\kappa}^{(2)}$ is $t$-acceptable, we have   for some $
v_{\max}^{(2)}$ and $ R_v^{(2)}\in \calP(g_1,\ldots, g_{t-1})$ with   all coefficients of size
$O(n^{o(1)}x)$ that
\begin{eqnarray}
 q'(y)
&=& \frac{d}{dy} \left(\sum_{v=0}^{\vmax^{(2)}}  R_v^{(2)}(\hat\zeta )
\left({y\over \la_t}\right)^{v}\right) -\frac{1}{y} \non \\
 &=&-{1\over y} + O\left(\frac{ n^{o(1)}x }{\lambda_t}\right)\lab{qkder}\\
&=&-{1\over j^*} + O\left(\frac{n^{o(1)}x}{\lambda_t}+\frac{|y-j^*|}{(j^*)^2}\right) \non
\end{eqnarray}
for $|y-j^*|=o(j^*)$,
 and on the other hand, from the definition of $j^*$, $q(j^*)=0$.
It   follows that for $k=j^*+ O(\sqrt {j^*} \log j^* )$, we have (again noting $j^*\sim \la_t$)
$$q(k)=\int_{j^*}^k q'(y)\, dy = -\frac{k-j^*}{j^*}  + o\left( (j^*)^{-1/2 } xn^{o(1)}\right).
$$
 
Thus, for the same range of $k$, summing~\eqn{Stexp} over $j$ between $\tilj:=\lfloor j^* \rfloor$
and $k$ gives \bel{ratioasy} \log( T_k/T_{\tilj})=\frac{-(k-\tilj)^2}{2j^*}+o(1) \ee (and this
argument applies whether $k$ is smaller or larger than $\tilj$). Hence, the sum of $T_k$ for $k =j^*+
O(\sqrt j^* \log j^* )$ is asymptotic to $T_{\tilj}$ times the sum of $e^{-(k-\tilj)^2/2j^*}$ over
the same range, and is hence
$$
( 2\pi j^*)^{1/2} T_{\tilj}(1+o(1)).
$$
Also,~\eqn{ratioasy} is valid at the extreme ends of the range, i.e.\ $k =j^*+ \Theta(\sqrt j^*
\log j^* )$. Thus, recalling~\eqn{ratio2}, sall the terms in~\eqn{Stm1d} outside the  range $k =j^*+ O(\sqrt j^* \log j^*
)$ are negligible and \bel{firstsum} \sum_{j=0}^{\lfloor  m_{t} \rfloor} T_j \sim ( 2\pi
\tilj)^{1/2} T_{\tilj}. \ee

To estimate $T_{\tilj}$, we use Stirling's formula and then $j^*\sim\la_t$ and $|\tilj-j^*|<1$ to
write \bel{Tjfirst} \frac{\la_t^{\tilj}}{\tilj!}\sim \frac{(e\la_t/\tilj)^\tilj}{\sqrt{2\pi\tilj}}
\sim  \frac{(e\la_t/j^*)^{j^*}} {\sqrt{2\pi\tilj}}. \ee Using~\eqn{jstarasymp} we may expand the
logarithm of $1/A_{t,\kappa}^{(3)}$ to obtain, for some acceptable polynomials $A_{t,\kappa}^{(4)}$ and
$A_{t,\kappa}^{(5)}$ in $\calP(\bfg_{t-1})$,
$$
\log(\la_t/j^*) =
 A_{t,\kappa}^{(4)} (\hat\zeta)-1+o(1/\la_t)
 $$
 and then
\bel{Tjsecond} 
(e\la_t/j^*)^{j^*}=\exp(j^*\log(e\la_t/j^*)) = \exp( \la_t A_{t,\kappa}^{(5)}
(\hat\zeta)+o(1)). 
\ee 
(Here $ A_{t,\kappa}^{(5)}$ just contains the significant terms of $ A_{t,\kappa}^{(3)}\cdot  A_{t,\kappa}^{(4)}$.)
Next, from~\eqn{xiexp} with $j=\tilj$ we have, for some $t$-acceptable polynomial $A_{t,\kappa}^{(6)}$,
\bel{A60}
\sum_{i=0}^{\tilj-1}\log\xi_{t,\eps}\left(f_{t,i}\right) = \la_t\big(
A_{t,\kappa}^{(6)}\big(\zeta(\tilj)\big)-1\big)+o(1).
\ee

For example, if $\xi_{t,\eps}$ happens not to contain $g_t$, then $A_{t,\kappa}^{(6)}$  is equal to  $1+g_t\, \widehat\log\,  \xi_{t,\eps}$, where $\widehat\log$ denotes the logarithm truncated to significant terms.
Since $|j^*-\tilj|<1$ and $ A_{t,\kappa}^{(6)}$ is $t$-acceptable, we may replace $\tilj$ in the right hand side of~ \eqn{A60}  by $j^*$,
with no other change to the equation. Using this, together with~\eqn{Tjfirst} and~\eqn{Tjsecond},
in~\eqn{Tdef} with $j=\tilj$, we may transform~\eqn{firstsum} into \bel{secondsum}
\sum_{j=0}^{\lfloor  m_{t} \rfloor} T_j \sim   \exp \Big( P_{t,\kappa} \big(\zeta(\tilj ) \big) +
\la_t  A_{t,\kappa}^{(6)}\big(\zeta( j^*)\big)-\la_t+\la_t A_{t,\kappa}^{(5)}\big(\hat\zeta\big)
\Big). 
\ee 
Note that $A_{t,\kappa}^{(6)} -1+A_{t,\kappa}^{(5)}$ is $t$-acceptable.   Then 
  the expansion~\eqn{jstarasymp} calls for   replacing $g_t$ in $A_{t,\kappa}^{(6)}$ by $A_{t,\kappa}^{(3)}(\hat\zeta )$:
$$
 A_{t,\kappa}^{(6)}\big(\zeta( j^*)\big)-1+A_{t,\kappa}^{(5)}\big(\hat\zeta\big)=
  A_{t,\kappa}^{(6)}\big(\zeta\big(\la_t A_{t,\kappa}^{(3)}(\hat\zeta ) \big)\big)-1+A_{t,\kappa}^{(5)}\big(\hat\zeta\big)=
 A_{t,\kappa}^{(7)}\big(\hat\zeta\big) + o(1/\la_t)
 $$
for some acceptable polynomial $A_{t,\kappa}^{(7)}\in \calP(\bfg_{t-1})$. Also, by hypothesis (ii) and the fact that
$|j^*-\tilj|<1$, we have $P_{t,\kappa} \big(\zeta(\tilj ) \big)=P_{t,\kappa} \big(\zeta(j^* )
\big)+o(1)$. Again replacing $g_t$ by $A_{t,\kappa}^{(3)}(\hat\zeta )$, using~\eqn{jstarasymp} we obtain
$$
P_{t,\kappa} \big(\zeta(\tilj ) \big)=\tilde P_{t,\kappa}( \hat \zeta  )+o(1)
$$
for a polynomial $\tilde P_{t,\kappa}\in \calP(\bfg_{t-1})$  that has exactly the same properties described in (ii) for
$ P_{t,\kappa}$. 

Note that there are multiple valid choices for $\tilde P_{t,\kappa} $ at this point, due to the possible inclusion of negligible terms. To avoid ambiguity, we specify that the  terms that are retained are exactly those that are significant in this argument when $p$ is precisely $n^{-\kappa}$, that is, terms 
 of order  $n^a p^b$ for which $a/b\ge k$.

Now from~\eqn{Stm1d} and~\eqn{secondsum} we have
\bel{Sfound}
S_{t-1}(j_1,j_2,\ldots,j_{t-1}) \sim  \exp \Big(\tilde P_{t,\kappa} ( \hat \zeta  ) +
\la_t A_{t,\kappa}^{(7)}\big(\hat\zeta \big)  \Big).
\ee
We may now set
$$
P_{t-1,\kappa} =  \tilde P_{t,\kappa} + \la_t A_{t,\kappa}^{(7)} 
$$
to obtain parts (i) and (ii) of the inductive hypothesis. Indeed, by this recursive definition we
obtain that
$$
P_{t,\kappa}=\sum_{t'=t+1}^s  \la_{t'}A_{t',\kappa}
$$
for some acceptable polynomials $A_{t',\kappa}$. Verifying part (iii) of the inductive hypothesis
requires simply noticing that the estimates in the above derivation are, inductively, uniform over
all appropriate $\hat\zeta$. This uses the uniformity of the estimates in Lemma~\ref{converge} and
Corollary~\ref{expanded}.

The inductive step is now fully established, and we have~\eqn{Sform} for all $t$. Taking
$t=0$,~\eqn{S0sim} shows that \bel{forkappa} \BP(X=0)\sim \exp( -P_{ 0,\kappa}). \ee By part
(ii) of the inductive hypothesis, $P_{0,\kappa}=(1+O(n^{o(1)}x))   \sum_{t\in\smalls} \la_{t}$.

We now show that \bel{equalpolys} \mbox{the polynomial $P_{0,\kappa}$ is a truncation of
$P_{0,\chi+\eps}$ for all 
$\chi+\eps\le\kappa < 2-\eps^{\prime\prime}$, } \ee 
(where the upper bound $2-\eps^{\prime\prime}$ arises from~\eqn{Egrows}). 
This statement immediately requires some
qualification. In the definition of $P_{t,\kappa}$, it is important to note that any expansions
during the proof above must be taken in the formal sense. For instance, if $\chi+\eps$ happens to
take certain rational values, then some terms in an expansion of the form $n^ap^b$ might  happen to be
equal to other terms $n^cp^d$, but these terms should be kept separate when comparing polynomials.

 We begin by showing that there is no ambiguity in the
definition of $P_{t,\kappa}$  due to the arbitrariness of ordering of the types in $\smalls$. That
is, we show that  the various orderings of types that are valid all lead to the same terms in
$P_{t,\kappa}$. Consider two possible orderings of types $\pi$ and $\tilde{\pi}$. For each choice
of ordering there corresponds a polynomial $P_{t,\kappa}$ in (\ref{forkappa}). Let us refer to the
function $o(1)$ in (\ref{pkappa}) as $g(n)$. Since $g(n)$ may be taken so that $n^{g(n)}$ is any
positive constant function, and for all such functions the two polynomials must have equal values
to within $o(1)$, all terms in the polynomials that are bounded below  when $n^{g(n)}$ is constant
must be equal. Terms that tend to 0 when $n^{g(n)}$ is constant must be $n^{-\eps'}$ for some
$\eps'>0$ and hence cannot occur in these polynomials.

We continue with the main part of the proof of~\eqn{equalpolys}. Note first, as an easy argument
shows, that as $\kappa$ increases smoothly from $\chi+\eps$ to $2-\eps''$,   there is a finite number of values of $\kappa$ at which the ordering of the
types can change, or a type changes from small to large.  (Recall that, as $\kappa$ increases, $p$
decreases, and hence every $\la_i$ decreases, and hence a type can move from $\smalls_1$ to
$\smalls_0$, and at essentially the same $\kappa$ from $\smalls_0$ to large, but not in the
reverse direction.) These are special values for our argument, since the ordering of types
determines the order of expansions in the inductive arguments concerning $S_t$. We designate the
minimum value, $\chi+\eps$, also as one of these special values, $\kappa_0$, and let the others be
$\kappa_1,  \kappa_2, \ldots$, with  $\kappa_0 < \kappa_1 < \cdots$.

Let us first fix two of these distinct values of $\kappa$,   $\kappa_i<\kappa_{i+1}$, and consider
$\kappa$ in the open interval $(\kappa_i,\kappa_{i+1})$. First, we will show that in the inductive
argument given above,  for such $\kappa$, we may use $P_{t,\kappa_i}$ in the argument in place of
$P_{t,\kappa}$ (subject to some near-trivial modification we will describe). We show moreover that
$P_{t,\kappa}$ is a truncation of $P_{t,\kappa_i}$. To be precise, we claim that all the
expansions in the argument for $\kappa$ can be replaced by the corresponding ones from the
argument for $\kappa_i$. The difference between the corresponding expansions lies only in the
terms that are absorbed by the error terms in the argument for $\kappa$. To see this inductively,
we need only to modify the argument for $\kappa$ slightly. We describe various aspects of the two
arguments as being ``for $\kappa$" or ``for $\kappa_i$" to distinguish between the two versions.

The inductive argument for $\kappa$ begins with a maximal $t\in\smalls$. Since $\kappa$ is not a
special value, it cannot be true that $t\in\smalls_0$. However, it may happen that a type $t'$ is
large for $\kappa$ but small (and hence in $\smalls_0$) for $\kappa_i$. By what has been shown
about ordering types arbitrarily, we may assume that types that are small for $\kappa$ have the
same ordering for $\kappa$ as they do for $\kappa_i$. For any type like the above-mentioned $t'$,
we may extend the definitions in the argument for $\kappa$ by putting $S_{t'}=1$, and it is easy
to verify that $P_{t',\kappa_i}=o(1)$ when evaluated at the value of $p$ occurring in the argument
for $\kappa$, i.e.\ $p=n^{-\kappa+o(1)}$. As the remaining types have identical order,  it remains
to be shown that if   $t\in \smalls$ for $\kappa$,  then $P_{t,\kappa}$ equals $P_{t,\kappa_i}$
except  for those terms of  $P_{t,\kappa_i}$ which are $o(1)$ for $\kappa$.

 At every point in the argument above for arbitrary $\kappa$ that an expansion is called for, beginning with the use of $\hatma$ in~\eqn{ratio}, we may add the extra terms called for in the $\kappa_i$ argument, and note that they fall into the error terms in the equation concerned. In particular, for~\eqn{ratio} this is true because of the assertion about the truncations  in Corollary~\ref{expanded}. Then, since this equation (and those following it) is true with these extra terms, the argument works as before, with expansions being carried out and with truncations determined by the argument   for $\kappa_i$ rather than $\kappa$. Every step of the argument then preserves the expansions obtained in the argument for $\kappa_i$, but all other aspects of the argument are as for $\kappa$. This is immediately obvious in places where products of series, and logarithms, are expanded, but it is a little more subtle in the part involving $j^*$, so we examine this in more detail.

We need to show that $A^{(3)}_{t,\kappa}$ equals $A^{(3)}_{t,\kappa_i}$ up
to insignificant terms. Let $\tilde{\la_t}$ be $\la_t$ with $p=n^{-\kappa_i+o(1)}$. Let
$\tilde{A}^{(3)}_{t,\kappa_i}$ be the polynomial derived with $p=n^{-\kappa_i+o(1)}$ but evaluated
at $p=n^{-\kappa+o(1)}$. We can write $A^{(3)}_{t,\kappa_i}$ as
$A^{(3)}_{t,\kappa_i}=B^{(3)}_{t,\kappa_i}+C^{(3)}_{t,\kappa_i}+D^{(3)}_{t,\kappa_i}$, where
$D^{(3)}_{t,\kappa_i}=o(1/\tilde{\la_t})$ and where $C^{(3)}_{t,\kappa_i}$ is significant for
$p=n^{-\kappa_i+o(1)}$ but such that $\tilde{C}^{(3)}_{t,\kappa_i}=o(1/\la_t)$. We constructed
$D^{(3)}_{t,\kappa_i}$ from a given number of contractions and the contraction constant is smaller
for $\kappa$ than it is for $\kappa_i$ (for the contractions obtained when the coefficients in
$A^{(2)}_{t,\kappa_i}$ and $A^{(2)}_{t,\kappa}$ are replaced by their absolute values), hence
$\tilde{D}^{(3)}_{t,\kappa_i}\leq D^{(3)}_{t,\kappa}$ and
$\tilde{D}^{(3)}_{t,\kappa_i}=o(1/\la_t)$. All of the remaining steps in the argument for
$p=n^{-\kappa+o(1)}$ involve sums, products, expansions of logarithms or substitutions into
polynomials and so everything arising from $\tilde{C}^{(3)}_{t,\kappa_i}$ is of the order
$o(1/\la_t)$. Thus, ignoring $o(1/\la_t)$ terms,
$\tilde{A}^{(3)}_{t,\kappa_i}=A^{(3)}_{t,\kappa}$.

Next, we will show that the inductive argument given above,  for   $\kappa_{i+1}$, remains valid
if we use $P_{t,\kappa}$ in the argument in place of $P_{\kappa_{i+1}}$, and   that
$P_{\kappa_{i+1}}$ is a truncation of $P_{t,\kappa}$. In this case, no type can move from being
small for the $\kappa$ argument to being large for the $\kappa_{i+1}$ argument (since $\kappa_{i+1}>\kappa$), but possibly a type $t$ is in $\smalls_1$ for the case of $\kappa$ and in
$\smalls_0$  for the case of $\kappa_{i+1}$. By part (ii) of the inductive hypothesis, the
contribution from the type $t$ to $P_{t,\kappa}$ is  $\la_t+o(1)$ when $p$ is taken in the
appropriate range for ${\kappa_{i+1}}$ because then $\la_t=n^{o(1)}$, and moreover this is also
the contribution to $P_{t,\kappa_{i+1}}$. The rest of the argument for this case only involves
considering the expansions, so is similar to the argument above.

Statement~\eqn{equalpolys} now follows by induction from the statements that  $P_{t,\kappa}$ is a truncation of $P_{t,\kappa_i}$ and   that $P_{\kappa_{i+1}}$ is a truncation of $P_{t,\kappa}$.
In view of the argument above that decreasing $\kappa$ simply adds more terms to  $P_{0,
\kappa}$, we see that decreasing $\eps$ does the same thing to $P_{ 0,\chi+\eps}$. Hence, this is
the truncation to a finite number of terms of a power series $F(G_0)$ in $n$ and $p$.
Since there is a bounded number of terms in~\eqn{maineq} that are $o(1)$ for a given $\kappa$, we have now established~\eqn{maineq} for this power series $F(G_0)$ and for $p=n^{-\kappa+o(1)}$ (whenever $\kappa\ge \chi+\eps$).  In particular, with the terms $c_\ell n^{i_\ell}p^{j_\ell}$   arranged in decreasing order of $i_\ell/j_\ell$, the claimed characterisation of $M_\eps$ follows.  Note that the function represented by $o(1)$ in (\ref{maineq}) is
given explicitly by  
$$
f(n,p)=\log\left(\BP(X=0)\right)- \sum_{\ell =0}^{M_\eps} c_\ell n^{i_\ell}p^{j_\ell}.
$$
We may now apply Lemma~\ref{panything} with
$a=\chi+\eps$ and $b=2-\eps^{\prime\prime}$
to deduce that the convergence in~\eqn{maineq} is uniform
over all $\kappa\in[\chi+\eps,2-\eps^{\prime\prime}]$.

All that remains is to show the  strict positivity of the exponents $i_\ell$ and $j_\ell$ in $F(G_0)$. Note that a term   $n^{i_\ell}p^{j_\ell} $ with $i_\ell\le 0$ must have $j_\ell<0$, otherwise it is always $o(1)$ and can simply be omitted. However, such a term is decreasing in $p$, so, if it is ever significant, must be so when  $p \le n^{-2+\eps''}$. However, at that point we know $\pr(X=0)\sim 1$, and hence the term must be insignificant here too. Thus, such terms can be dropped. It follows that we may assume $i_\ell>0$. Given (by the same argument) that the term must be insignificant for small $p$, we deduce that $j_\ell>0$ also.  
The  $\Gnp$ case of the theorem follows.
  \qed

\section{Graphs with forbidden subgraphs in ${\cal G}(n,m)$}\lab{s:maingnm}

We will show that the $\Gnp$ case of Theorem~\ref{t:main} can be extended to give a similar
result in $\Gnm$ without much difficulty. Specifically, we provide
asymptotics for the probability of $\Gnm$ not containing a fixed subgraph
isomorphic to $G_0$.
The asymptotics could be expressed in terms of $n$ and $m$, but
it is more convenient to use $n$ and  the parameter
$d =m / {n \choose 2}$  defined in
(\ref{ddef}).  We employ the $\Gnp$ case inside the proof, for a value of $p$ that is close, but not quite equal, to $d$, though for the statement of the theorem we have renamed $d$ as $p$ for convenience.
\bigskip

 \noindent
{\bf Proof of the $\Gnm$ case of Theorem~\ref{t:main}.}

Let $Y$ denote the number of edges of a graph. The probability that $X=0$ in $\Gnm$ is precisely  $\pr(X=0\mid Y=m)$ in $\Gnp$. In the rest of the proof we estimate this quantity, with all probabilities referring to $\Gnp$. By Bayes' Theorem, what we desire is
\bel{rewrite}
\pr(X=0\mid Y=m) = \pr(Y=m\mid X=0)\frac{\pr(X=0)}{\pr(Y=m)}
\ee
This formula  is valid for all $0<p<1$. The value of $p$ we will use, which is  specified below, is asymptotic to $d$ and hence lies in the range required for the 
$\Gnp$ case of Theorem~\ref{t:main}, given by~\eqn{pkappa} with the same restrictions on $\kappa$, which determines $\smalls$ via~\eqn{smallsdef}. 
Thus, Theorem~\ref{t:main} gives us $\pr(X=0)$ in $\Gnp$. 

The main difficulty is
computing $\pr(Y=m\mid X=0)$. For this,
 we will first alter the analysis in Section~\ref{s:main} to consider  the $G_0^*$-clustering in  $\Gnp$. Recall that this
  is obtained by adding  to
$\smalls$ the type $t^*$ of maximal cluster corresponding to a single edge. For convenience, we henceforth denote the cluster type $t^*$ by 0.  

Considering the polynomial $\xi_{0,\epsilon}(n,p,\bfg)$ in Corollary~\ref{expanded}, for
$j/\lambda_0\leq 3$ (in accordance with (\ref{fbound})), by part (c) of that Corollary  
\bel{hatma0approx}
\hatma_0(j\delta_0)=\xi_{0,\epsilon}(n,p,\hat\bfg(j))+o(\lambda_0^{-1}),
\ee
where  $\hat g_0(j)=j/\lambda_0$,  $\hat g_i(j)=0$ for $i\geq 1$,
 provided that $p=p(n)=O(n^{-\chi-\epsilon})$ and satisfies \eqn{Egrows}.

 Also define   $\tilde\bfg$ by $\tilde g_0=m/\la_0$ and $\tilde g_i=0$ for $i\geq 1$, and let $\xi $ denote $\xi_{0,\epsilon}(n,p,\tilde\bfg)$ (noting that $\xi$ is a function of $n$, $p$  and $m$.
 As $\hat g_0(j)\le 3=O(1)$ we have 
\bel{xiapprox}
\xi_{0,\epsilon}  \big(n,p,\hat \bfg(j)\big)= \xi +O\big( x(m-j)/\la_0\big) 
\ee
by Corollary~\ref{expanded}(a).   

By the definitions of $\rho(f,h)$ and $\gamma(f,t)$
before Proposition~\ref{thetabound2}, 
one would expect that the probability that $\Gnp$ has no copies of $G_0$
and $m'$ edges will
be maximised, given $p$,  at $m'\approx m$  provided that $\rho(m \delta_0,\delta_0)\approx 1$, or 
 $\gamma(m\delta_0,0)\approx m/\lambda_0$. On the other hand, in $\Gnp$ the ratio of  the probabilities of having a given number of edges, when increasing $m$ to $m+1$, is approximately $d/p$. Consequently,
we define $p$ by
\bel{mudef2}
 p =d/ \xi 
 \ee
 (recalling that~\eqn{ddef} gives
$d$ as a function of $n$ and $m$).
Then~\eqn{hatma0approx} and Lemma~\ref{converge}   imply that 
\bel{pd}
 p=d\big(1+O(x+\la_0^{-1})\big),
 \ee
 and hence our assumptions on $d$ imply the necessary properties of $p$ such as~\eqn{Egrows}, perhaps with different values of the unimportant constants. 
 
From the $\Gnp$ case of Theorem~\ref{t:main},  $\pr(X=0)$ in $\Gnp$ is $e^{-\Theta(\la_t)}$, and $\la_t=o(\la_0)$  by~\eqn{meanratio}.  
On the other hand,  
The number $Y$ of edges  in ${\cal G}(n,p)$ is distributed as ${\rm Bin} (N,p)$ where $N=\binom{n}{2}$,
with mean $\lambda_0=Np\sim m$. Hence, $\pr(Y>2m)<e^{-cm}<e^{-\Omega(\la_0)}$ (for instance by Chernoff's bound).
It follows that
\bel{mainsum}
\pr(X=0)\sim\sum_{j\le 2m} \BP(\class_{j\delta_0}).
\ee

Using the definition of $\rho$, Proposition~\ref{hatmaclose} and Lemma~\ref{converge}, and then~\eqn{hatma0approx} and~\eqn{xiapprox}, we have  
\bea
\rho(j\delta_0,\delta_0)&=&\frac{\la_0}{j+1}\gamma (j\delta_0,0)\nonumber\\
&=&\frac{\la_0\hatma_0(j\delta_0)}{j+1}\big(1+o(\la_0^{-1})\big)\nonumber\\
&=&\frac{\la_0\xi }{j+1}\big(1+ o(\la_0^{-1})+O( x(m-j)/\la_0) \big)\nonumber\\
&=&\frac{m}{j+1}\exp\big( o(\la_0^{-1})+O( x(m-j)/\la_0) \big) \nonumber
\eea 
by~\eqn{mudef2} and~\eqn{ddef}.
Hence~\eqn{mainsum} gives
\begin{eqnarray*}
 \frac{\BP(X=0)}{\BP(Y=m,X=0)}&=&
\frac{\BP(X=0)}{\BP(\class_{m\delta_0})}\\
&\sim&
\sum_{j\le 2m}
\frac{\BP(\class_{j\delta_0})}
{\BP(\class_{m\delta_0})}\\
&=&
\sum_{j\le 2m}
\rho(m\delta_0,(j-m)\delta_0)\\
&=&
\sum_{j=  m}^{2m}\prod_{i=m}^{j-1}
\rho(i\delta_0,\delta_0)
+\sum_{j=  0}^{m}\prod_{i=j}^{m-1}
\rho(i\delta_0,\delta_0)^{-1}
\\
&=&
\frac{m!}{m^m\ }
\sum_{j=0}^{2m}
\frac{m^j }{j!}\exp\big(o((m-j)/\la_0)+O(x(m-j)^2/\la_0)\big)\\
&\sim&
\frac{m!}{m^m\ }
\sum_{j=0}^{2m} 
\frac{m^j }{j!}\\
&\sim&
\frac{m!e^{m }}{m^m }\sim 
\sqrt{2\pi m}.
\end{eqnarray*}
In the third-last line, the main terms of the summation have $|m-j|\approx \sqrt m\sim \sqrt  {\la_0}$, for which the error terms are $o(1)$ as $x\to 0$. The remaining terms are insignificant since the absolute value of the $j$th term in the sum is $ \frac{m^m}{m!}\exp(-\Omega(m-j)^2/\la_0)$, 
which dominates the error term. 
 The last line uses Stirling's formula.

Taking the multiplicative inverse of the previous asymptotic formula produces 
\eq\lab{X0kcond2}
\BP(Y=m   \mid  X=0)\sim\frac{1}{\sqrt{2\pi m}}.
\en

For the other  factors in~\eqn{rewrite}, first recall that $\xi$ comes ultimately as a truncation of the power series $\xi_0$, in $n$ and $p$ (here $t=0)$ in Corollary~\ref{expanded}. Thus, we can use~\eqn{mudef2} and~\eqn{pd} to expand $p$ as a power series in $n$ and $d$. Specifically, we obtain $p=d\tilde J_1\big(1+o(\la_0^{-1})\big)$  where $\tilde J_1$ is the truncation of a power series $J_1$ in $n$ and $d$ to significant terms. Here $J_1$ is independent of $\kappa$, being the termwise limit of the power series obtained for $\kappa$  as $\kappa\downarrow\chi$ (which represents increasing $p$).  This can be  substituted into the polynomial obtained by truncating
the  power series for  $\log\pr(X=0)$ 
  obtained from the $\Gnp$ case of Theorem~\ref{t:main},
at an appropriate level,
to express $\log \pr(X=0)$ as $ \tilde J_2+o(1) $ where $\tilde J_2$ is 
 a truncation of a power series $ J_2$ in $n$ and $d$, with $J_2$ independent of $\kappa$.      Similarly, $\pr(Y=m)$ is simply the binomial probability which can be estimated using Stirling's formula. The leading (polynomial-type) factor  is asymptotic to  $1/\sqrt {2\pi m}$, which cancels with $\BP(Y=m|X=0)$ obtained above. The logarithm of the exponential factor can be expanded using $p=d\tilde J_1\big(1+o(\la_0^{-1})\big)$   to obtain an expansion of the type required. (Theorem~\ref{t:Gnm} gives an example.) Subtracting this from the expansion for $ \log\pr(X=0)$  gives the $\Gnm$ case of the theorem by~\eqn{rewrite}. The positivity of the 
 exponents $i_\ell$ and $j_\ell$  follows by arguing as in the proof of the $\Gnp$ case. \qed

\newpage

\begin{appendices}
 \section{ Calculations for triangle-free graphs}\lab{s:tri}

\noindent{\bf Proof of Theorem~\ref{t:Gnp}}

Section~\ref{s:main} shows that an asymptotic formula for the probability
a subgraph $G_0$ is not present in $\Gnp$ exists, but it does not state
the formula explicitly. Nevertheless, the proof  fully prescribes a method of
calculating the formula for any particular case.  At its heart, the proof uses Corollary~\ref{expanded}, in which the power series $\xi_{t,\eps}$ are not stated explicitly. To obtain a formula in practice, these must be determined
to a required accuracy, along with the quantities $c(u,t,h)$  defined in (\ref{cdef}).
In this section
we demonstrate how the necessary calculations are performed in the case when
$G_0$ is a triangle.

\begin{figure}[h]
\begin{center}
\includegraphics[width=14cm,trim=0cm 0cm 2.2cm 12.5cm,clip=true]{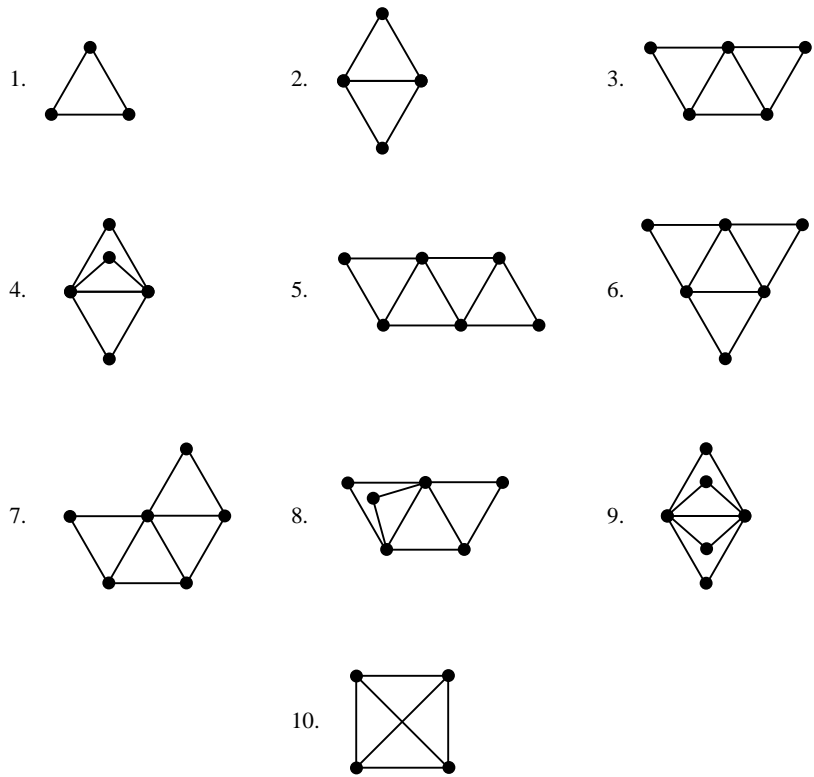}
\end{center}
\caption{\it Ten types of cluster
}
\label{f:types}
\end{figure}

We let $G_0=K_3$, the complete graph on 3 vertices, and proceed to estimate the probability that
${\cal G}(n,p)$ contains no triangles in the case that $p<n^{-7/11-\eps}$. (This constraint will be relaxed to $p=o(n^{-7/11})$ at the end.) It is easy to check that
that \eqn{chidef} determines $\chi = \frac12$ when $G_0=K_3$.
If we make the restriction
$p=n^{-\kappa+o(1)}$ with $\kappa> 7/11$,  then there are then 10 possible cluster types possible in $\smalls$ according to~\eqn{smallsdef}. We thus have    $\smalls=\{1,2,\ldots,10\}$ as depicted in
Figure~\ref{f:types}.

All
these types are present in $\smalls$ when $\kappa$ is
at most 2/3.  All  other cluster types have expected number tending to
$0$ as $\kappa> 7/11$, and are
therefore are not in
$\smalls$.   Recall that
the poset ordering $\prec$ on $\smalls$ is not necessarily a linear ordering; for example,
the types $\{5,6,7,8,9,10\}$ are all maximal,  and therefore not comparable.
The ordering $\prec$ is extended to the usual real linear ordering on $\smalls$
denoted by $<$.

The first step  is to calculate $\la_t$ for $t\in\smalls$.
In accordance with~\eqn{tsize} and~\eqn{lambdadef}, we obtain the $\la_t$ as in Table~\ref{tab:lambdas}.
\begin{table}[h]
\vspace{1cm}

\begin{center}
\begin{tabular}{|c|c|c|c|c|}
\hline
\raisebox{-0.4cm}{\rule{0cm}{10mm}}$\la_1=\frac{1}{6}[n]_3p^3$&
$\la_2=\frac{1}{4}[n]_4p^5$&
$\la_3=\frac{1}{2}[n]_5p^7$&
$\la_4=\frac{1}{12}[n]_5p^7$&
$\la_5=\frac{1}{2}[n]_6p^9$\\
\hline
\raisebox{-0.4cm}{\rule{0cm}{10mm}}$\la_6=\frac{1}{6}[n]_6p^9$&
$\la_7=\frac{1}{2}[n]_6p^9$&
$\la_8=\frac{1}{2}[n]_6p^9$&
$\la_9=\frac{1}{48}[n]_6p^9$&
$\la_{10}=\frac{1}{24}[n]_4p^6$\\
\hline
\end{tabular}
\end{center}
\caption{Expected numbers of small clusters.}
\lab{tab:lambdas}
\end{table}

Our next task is to find the polynomial $\xi_{t,\eps}$ of Corollary~\ref{expanded} for all $t\in\smalls$. For this, the proof of the corollary describes an iterative scheme to compute the $F_t^{(r)}$ and hence $\hatma^{(r)}_t$.

We can drop all terms that would yield coefficients of variables $g_{t_i}$ that are $O(\error/\lambda_t)$  for some $\epshat>0$. This is because in $\xi_{t,\eps}(f)$, each $g_{t_i}$ is assigned a  value that is  $n^{o(1)}$, and hence the dropped terms are subsumed into the error term in~\eqn{excor} when $\epshat$ is sufficiently small (recalling that $\hatma_t(f)\sim 1$ by Lemma~\ref{converge}).
For similar reasons, we can drop any $O(p^2)$ term in the expansion of $c(t,t,\bz)$ at the front of~\eqn{Ftdef}, as  $p^2\la_t=O(p^2\la_1) = O(p^5n^{3})=o(n^{-2/11})$. Note that
$ c(1,1,\bz)=1-p^3$
since it is simply the probability that three vertices do not form a triangle. Hence, we can treat $c(1,1,\bz)$ as 1. A similar argument applies to
  $
c(t,t,\bz)
$
for all other $t\in\smalls$.

Moving on to the quantities $c(u,t,h)$ inside the summation in~\eqn{Ftdef},   for any non-zero $  h\in\functions$, clearly
$c(t,t,h)=O(p^3)$,
and so these terms can  be ignored completely for the same reason, for all $t$.

For the other terms in the summation, we only need to compute $ c(u,t,h)$ to $O(\error/\la_u)$. For $u=2$, note that $ \error/n^4p^5 =\Omega(p^{10/7})$ for sufficiently small $\epshat$ since $n<p^{-11/7}$. Thus, we may drop $p^2$ terms in $c(2,t,h)$. First
consider  $c(2,1,\bz)$. In~\eqn{cdef}, $J$ is a cluster of type 2, i.e.\ (the edge set of) two triangles with a common edge. $Q$ corresponds to one of the two triangles of $J$ (so there are two choices for $Q$).  There are four cases for $H$, as it must contain $J\setminus Q$ but no triangles. Letting $q=1-p$, we get
\[
c(2,1,\bz) = 2(q^3+2pq^2+p^2q) = 2(1-p)+O(p^2).
\]
The other cases of $c(2,t,h)$ can be computed similarly, and only $h=\delta_1$ is significant (i.e.\ not $O(p^2)$).
Similarly, for $u=3$, $1/n^5p^7=\Omega(p^{6/7})$ and we may drop the $O(p)$ terms. The same clearly holds for all $u>3$ as well. In this way, we obtain all significant terms of $c(u,t,h)$  for $u>t$, as shown in   Table~\ref{tab:cs}.
In computing these, note that $h$ is quite restrictive. For instance, for $c(3,1,\bz)$, the   deletion of $Q$ from $J$ must  leave  no triangles, and there is only one such choice for $Q$.

\def\fracspace{\raisebox{-2mm}{\rule{0cm}{6mm}}}
 \def\exspace{\fracspace}
\begin{table}[h]
 \small
\begin{center}
\begin{tabular}{|c|c|c|c|c|}
\hline
\exspace  $u$& t&h&$c(u,t,h)$& cofactor\\
\hline
\fracspace 2&1&$\bz$& $ 2(1-p) $& $\frac 32 np^2 \cdot \hatma_2 $\\
\fracspace 2&1&$\delta_1$& $ 2p $& $\frac 32 np^2 \cdot  \hatma_2 g_1 \hatma_1^{-1}  $\\
\hline
\exspace 3&1&$\bz$& $ 1 $& $3n^2p^4 \cdot  \hatma_3 $\\
\exspace 3&1&$\delta_1$& $ 2 $& $3n^2p^4  \cdot  \hatma_3 g_1 \hatma_1^{-1}  $\\
\hline
\fracspace  4&1&$\bz$& $ 3 $& $\frac12 n^2p^4 \cdot   \hatma_4 $\\
\hline
\exspace 5&1&$\delta_1$& $ 2 $& $3n^3p^6 \cdot  g_1 $\\
\exspace 5&1&$\delta_2$&$2$    &$3n^3p^6  \cdot  g_2$ \\
\hline
\exspace 6&1&$\bz$      & $1$ & $n^3p^6  $\\
\exspace 6&1&$2\delta_1$& $3$ & $n^3p^6  \cdot  g_1^2$\\
\hline
\exspace7&1&$\delta_1$& $2$ & $3n^3p^6  \cdot  g_1   $\\
\exspace7&1&$\delta_2$& $2$ & $3n^3p^6  \cdot  g_2 $ \\
\hline
\exspace8&1&$\bz$& $1$& $3n^3p^6   $\\
\exspace 8&1&$\delta_1$& $2$ & $3n^3p^6  \cdot  g_1  $\\
\exspace 8&1&$\delta_2$& $1$ & $3n^3p^6   \cdot g_2  $\\
\hline
\fracspace 9&1&$\bz$& $4$& $\frac{1}{8}n^3p^6  $\\
\hline
\fracspace 10&1&$\bz$& $4$& $\frac{1}{4} np^3  $\\
\hline
\end{tabular}
\quad
\begin{tabular}{|c|c|c|c|c|}
\hline
\exspace 3&2&$\bz$&    $2$&  $2np^2\cdot  \hatma_3 $\\
\hline
\fracspace 4&2&$\bz$&   $ 3 $&  $\frac13 np^2\cdot  \hatma_4$\\
 \hline
\exspace  5&2&$\bz$      & $1$& $2n^2p^4   $\\
\exspace 5 &2&$\delta_1$ & $2$&  $2n^2p^4\cdot  g_1 $ \\
 \hline
\fracspace  6&2&$\bz$      & $3$ & $\frac23 n^2p^4  $\\
 \hline
\exspace  7&2&$\bz$      & $1$& $2n^2p^4  $\\
\exspace 7&2&$\delta_1$  & $2$&  $2n^2p^4\cdot  g_1  $ \\
\hline
\exspace 8&2&$\bz$      & $3$ &  $2n^2p^4    $\\
\exspace 8&2&$\delta_1$ & $1$ &    $2n^2p^4\cdot   g_1 $ \\
\hline
\fracspace 9&2&$\bz$      & $6$&   $\frac{1}{12}n^2p^2  $\\
\hline
\fracspace  10&2&$\bz$     & $6$&  $\frac{1}{6}p  $\\
\hline
\exspace 5& 3 &$\bz$   &   $2$& $ np^2  $\\
\hline
\fracspace 6&3 &$\bz$ &      $3$& $\frac13 np^2  $\\
\hline
\exspace 7&3 &$\bz$ &     $2$& $ np^2  $\\
\hline
\exspace 8 & 3 &$\bz$     & $2$&    $ np^2 $ \\
\hline
\exspace 8&   4 &$\bz$    & $1 $ & $6np^2  $\\
\hline
\fracspace 9&   4 &$\bz$     & $4$ & $\frac14 np^2  $\\
\hline
\end{tabular}

\end{center}
\caption{Significant contributions to~\eqn{Ftdef}}
\lab{tab:cs}
\end{table}

The   ``cofactor'' column of   Table~\ref{tab:cs} shows the significant contribution to those terms in $F_t$ from
$$
\frac{\lambda_u}{\lambda_t} \hatma_u
\prod_{i=1}^{k} \frac{g_{t_i}}{\hatma_{t_i}}.
$$
Here, and in the rest of the calculation, we assume $\epshat>0$ is as small as we like, and any  terms that are $O(\error/\lambda_t)$   are dropped. In each case, $\lambda_u /\lambda_t$ is the first item in the column, with any others (that are not equal to 1) appearing after ``$\cdot$''.  In each case only the leading term of $\lambda_u /\lambda_t$ turns out to be significant, since the correction terms are $O(1/n)$ and $\la_u/n=O(n^3p^5)$ for $u\ge 2$.  Any other factors which appear to be missing have simply been replaced by 1, with the following justification. In the initial iteration, for computing $F_t^{(1)}$ we have all $\hatma_v$ equal to 1, and by induction, thereafter they are $1+O(np^2)$ (if we treat each $g_{t_i}$ as 1). In the end each $g_{t_i}$ is substituted by something that is $n^{o(1)}$. Hence we may set any $\hatma_u$ or $\hatma_{t_i}$ equal to 1 in all iterations for all $u\ge 5$, since then  $\la_u np^2= O(n^7p^{11}+ n^5p^8) =O(\error/\lambda_t)$. Of course there are no contributions from $t\ge 5$ since all such $t$ are maximal in $\smalls$, and  $c(u,t,h)=0$ unless  $t\prec u$ (and we have already dealt with the case $u=t$).

The significant terms of~\eqn{Ftdef} are now deduced to be
\begin{eqnarray*}
F_1&=& - np^2
\left(3(1-p)\hatma_2+3p\hatma_2g_1/\hatma_1 \right) -  n^2p^4 \left(3\hatma_3+6\hatma_3g_1 /\hatma_1 +\frac{3}{2} \hatma_4\right)\\
 &&
-  n^3p^6\left(  18g_1+15g_2 +3g_1^2  +\frac{9}{2}\right) - np^3,\\
F_2&=& -
np^2\left(4\hatma_3+\hatma_4\right)- n^2p^4\left( 10g_1+25/2\right)-p,\\
F_3&=& -7np^2,\\
F_4&=& -7np^2,\\
F_t&=& 0  \quad \mbox{($t\ge 5$)}.
\end{eqnarray*}
Write $y=np^2$ and solve~\eqn{Ftdef} iteratively as described after that equation. It may help to note that any terms   of order $yp^2$, $y^2p$ or $y^4$ can   be dropped. After three iterations (actually the expressions don't change after the second update), the error is of order $x^4=\max\{y^4,p^4\}$ by~\eqn{x}, which is neglibible for each $t$. This gives $\xi_{t,\eps}  =1+F_t^{(4)}$ given as follows. \begin{eqnarray}
\xi_{1,\eps} &=&
1-3y+5py-3g_1py+\frac{21}{2}y^2
-6g_1y^2-\frac{81}{2}y^3+36g_1y^3-3g_1^2y^3-15g_2y^3,\lab{hatma1}\\
\xi_{2,\eps}&=&
1-p-5y+\frac{45}{2}y^2-10g_1y^2,\lab{hatma2}\\
\xi_{3,\eps}&=&1-7y,\lab{hatma3}\\
\xi_{4,\eps}&=&1-7y.\lab{hatma4}\\
\xi_{t,\eps}&=&1 \quad\mbox{ ($t\ge 5$)}.\lab{maxsim}
\end{eqnarray}

 We will evaluate the expressions given in Section~\ref{s:main}
with $7/11<\kappa<2/3$, so that
$\smalls_1=[10]$ and $\smalls_0=\emptyset$  (and actually $y=x$ as per~\eqn{xdef}). The ultimate result will then be valid for all values of $\kappa > 7/11$ by~\eqn{equalpolys}. We also fix $\eps$ in the range
$0<\eps<7/11-\chi=3/22$. With $\kappa$ and $\eps$ in these ranges, the
$\xi_{t,\eps}$ are given by the expressions 
(\ref{hatma1})--(\ref{maxsim}).

The recursive definition~\eqn{Sdef} of $S_t$ for $t\le 10$ starts with
  $S_{10}=1$ and hence, in~\eqn{Sform}, $P_{10,\kappa}=0$. Hence (just before~\eqn{ratio2}) $A_{10,\kappa}^{(1)}=1$. Of course there are options in choosing $A$'s since they are only determined up to an error term; we use the natural choices.

The next step is to determine $S_9$ and $P_{9,\kappa}$.
From~\eqn{A2def}, $A_{10,\kappa}^{(2)}=1$. Now (\ref{jstardef}) implies $j^\ast=\la_{10}$ 
and hence from~\eqn{jstarasymp}
  $A_{10,\kappa}^{(3)}=1$.
It is now easy to check that $A_{10,\kappa}^{(4)} =A_{10,\kappa}^{(5)}=1$ at~\eqn{Tjsecond},
and then similarly  $A_{10,\kappa}^{(6)}=A_{10,\kappa}^{(7)}=1$. (Much more detail in the steps here is provided in the less trivial
case when $t=1$ below.)
Finally, we conclude that, at~\eqn{Sfound},
$S_9(j_1,j_2,\ldots,j_9)\sim e^{\la_{10}}$
and then $P_{9,\kappa}=\la_{10}$.
In the same way one can show that
$S_t \sim \exp\left(\sum_{u=t+1}^{10}\la_u\right)$  for $t=8,7,6,5,4$.
In particular we have
$S_4\sim\exp\left(\sum_{u=5}^{10}\la_u\right)$ and
$P_{4,\kappa}=\sum_{u=5}^{10}\la_u$.

Next consider  $S_3$ and $P_{3,\kappa}$.
We have that $P_{4,\kappa}(\zeta)$ is independent of $\zeta$, and so $A_{4,\kappa}^{(1)}=1$. The ratio in (\ref{ratio}) is
$
T_j/T_{j-1}=(1-7y)\la_4(1+O(\eta))/j,
$
so 
$A_{4,\kappa}^{(2)}=1-7y$, 
$j^\ast=(1-7y)\la_4$ and  $A_{4,\kappa}^{(3)}=1-7y$.
Moreover, $\log(\la_4/j^\ast)=-\log(1-7y)=7y+O(y^2)=7y+O(\la_4^{-1})$, so
$A_{4,\kappa}^{(4)}=7y$ and
$$
\left(\frac{e\la_4}{j^\ast}\right)^{j\ast}=
\left(\frac{e}{1-7y}\right)^{\la_4(1-7y)}
=e^{\la_4+o(1)}
$$
from which we deduce $A_{4,\kappa}^{(5)}=1$.
Now,
$\sum_{i=0}^{\tilde{j}-1}\log \xi_{4,\eps}=\tilde{j}\log(1-7y)
=-7y(1-7y)\la_4+o(1)=-7y\la_4+o(1)$ implies that
$A_{4,\kappa}^{(6)}=  1-7y g_4$ and $A_{4,\kappa}^{(7)}=1-7y$.
Because ${P}_{4,\kappa}$ does not depend on $g_4$,
$\tilde{P}_{4,\kappa}=P_{4,\kappa}=\sum_{u=5}^{10}\la_u$. Finally, we have
$$S_3(j_1,j_2,j_3)\sim\exp\left(\sum_{u=5}^{10}\la_u+(1-7y)\la_4\right)$$
and $P_{3,\kappa}=\tilde{P}_{4,\kappa}+\la_4A_{4,\kappa}^{(7)}
=\sum_{u=5}^{10}\la_u+(1-7y)\la_4$.
Similar analyses which we   omit
show that 
$$
S_2(j_1,j_2) \sim \exp\left((1-7y)(\la_3 +\la_4)+  \sum_{u=5}^{10}\la_u \right).
$$
  Next, note that   $A_{2,\kappa}^{(1)}=1$,  $A_{2,\kappa}^{(2)}=A_{2,\kappa}^{(3)} = \xi_{2,\eps}$, $j^*=\la_2 \xi_{2,\eps}$,  $A_{2,\kappa}^{(4)}$ is a truncation of the expansion of $1-\log \xi_{2,\eps}$, 
$A_{2,\kappa}^{(5)}=1-\frac{25}{2}y^2$, $A_{2,\kappa}^{(6)}$ is the truncation of $1+g_2\log \xi_{2,\eps}$, which is 
$1+g_2 (-p-5 y+(10-10 g_1) y^2)$, and  $A_{2,\kappa}^{(7)}=1-p-5y+\frac{45}{2}y^2-10g_1y^2$.
  Eventually
$S_1(j_1)\sim\exp\left(P_{1,\kappa}(\zeta(j_1))\right)$
 where
\eq\lab{P1}
P_{1,\kappa}=
\left(1-p-5y+\frac{45}{2}y^2-10g_1y^2\right)\la_2 +(1-7y)(\la_3 +\la_4) +  \sum_{u=5}^{10}\la_u.
\ee

The final step of the induction is a little more involved.
We have 
\begin{eqnarray*}
\exp\big(P_{1,\kappa}(\zeta(j)) -
 P_{1,\kappa}(\zeta(j-1))\big)
&=&\exp\left(\frac{-10y^2\la_2}{\la_1}\right)\\
&=&\exp\left(  -15y^3+O(y^3/n)\right)
\end{eqnarray*}
and hence    $A^{(1)}_{1,\kappa}=1-15y^3$.
For~\eqn{A2def}  we set  $g_2=0 $ 
 to get $\tilde \xi_{1,\epsilon}$ and obtain
 $$
A_{1,\kappa}^{(2)}=1+c_1+c_2g_1
+c_3g_1^2.
$$
where
\begin{eqnarray*} c_1&=&-3y+5py+\frac{21}{2}y^2-\frac{81}{2}y^3 -15y^3\\
&=& -3y+5py+\frac{21}{2}y^2-\frac{111}{2}y^3,\\
c_2&=&-3py-6y^2+36y^3\\
c_3&=&-3y^3,
\end{eqnarray*}

The equation~(\ref{jstardef}) for $j^*$  becomes
\bel{iter}
{j^*\over\la_1}=1+c_1+c_2\left(\frac{j^*}{\la_1}\right)+
c_3\left(\frac{j^*}{\la_1}\right)^2.
\ee
Since the $c_i$ are $O(y^i)$, $j^*\sim \la_1$ and $\la_1y^4=o(1)$, we find iteratively that $
j^*=(1+c_1+c_2+c_1c_2+c_3)\la_1+o(1)$, and so
\eq\lab{A3}
A_{1,\kappa}^{(3)}=1+c_1+c_2+c_1c_2+c_3.
\en
Expanding $1-\log A_{1,\kappa}^{(3)}$ gives
$$
A_{1,\kappa}^{(4)}=1-c_1-c_2-c_3+c_1^2/2-c_1^3/3 
$$
and then truncating  $A_{1,\kappa}^{(3)}\cdot A_{1,\kappa}^{(4)}$ 
  gives
$$
A_{1,\kappa}^{(5)}=1-c_1c_2-c_1^2/2+c_1^3/6.
$$
\ Next, referring to~\eqn{A6} and writing $\tilde c_1= c_1 +15y^3 $, 
\begin{eqnarray*}
\sum_{i=0}^{\tilde{j}-1}
\log\xi_{1,\epsilon}(f_{1,i})
&=&
\sum_{i=0}^{\tilde{j}-1}
\log\left(1+\tilde c_1  +  {c_2 i\over\la_1}  +  {c_3i^2\over\la_1 ^2}\right)
 \\
&=& o(1)+
\sum_{i=0}^{\tilde{j}-1}
\bigg( \tilde c_1-\frac {\tilde c_1^2}{2} + \frac{\tilde c_1^3}{3}+\frac{(c_2   -\tilde c_1 c_2) i}{\la_1}+\frac{c_3 i^2}{\la_1^2 } \bigg)\\
&=&
\left[
\left(\tilde c_1-{1\over 2} \tilde c_1^2 +{1\over 3} \tilde c_1^3 \right){\tilde{j}\over\la_1}
+\left({1\over 2}c_2-{1\over 2} \tilde c_1c_2\right)
\left({\tilde{j}\over\la_1}\right)^2+
{1\over 3}c_3\left({\tilde{j}\over\la_1}\right)^3
\right]\la_1+o(1),
\end{eqnarray*}
so
\eq\lab{A6}
A_{1,\kappa}^{(6)}=
1+\left(c_1-{1\over 2}c_1^2 +{1\over 3}c_1^3 \right)g_1
+\left({1\over 2}c_2-{1\over 2}c_1c_2\right)
g_1^2
+{1\over 3}c_3 g_1^3.
\en

Substituting (\ref{A3}) for $g_1$ in (\ref{A6}),  dropping insignificant terms, and adding $-1+A_{t,\kappa}^{(5)}$,  we obtain  after some algebra 
 \begin{eqnarray*}
A_{1,\kappa}^{(7)}&=&
 1+\tilde c_1 +\tilde c_1 c_1-\frac12 \tilde c_1 ^2+\frac12 c_2+\frac12 \tilde c_1  c_2+\frac13 \tilde c_1 ^3-\frac12 \tilde c_1 ^2 c_1+\frac13 c_3- \frac12 c_1^2+\frac16 c_1^3+o(\la_1^{-1})\\
&=&
1-3y+\frac{7}{2}py+\frac{15}{2}y^2-\frac{29}{2}y^3.
\end{eqnarray*}
 
 Since $j^*=\la_1(1+O(y))$, changing  from $g_1= \tilde{j}/\la_1$ or $g_1= {j}^*/\la_1$ to $g_1=1$  in (\ref{P1})
induces a change to $P_{1,\kappa}$ of order
$O((j^\ast-\la_1)y^2\la_2/\la_1)=O(y^3\la_2)=o(1)$
and therefore $\tilde{P}_{1,\kappa}=P_{1,\kappa} |_{g_1=1}$.
Finally,
\begin{eqnarray*}
P_{ 0,\kappa}&=&
\tilde{P}_{1,\kappa}(1)+\la_1A_{1,\kappa}^{(7)}\\
&=&
\sum_{u=5}^{10}\lambda_u+(1-7y)\la_4+(1-7y)\la_3+
\left[1-p-5y+\frac{25}{2}y^2\right]\la_2\\
&&+\,
\left[
1-3y+\frac{7}{2}py+\frac{15}{2}y^2-\frac{29}{2}y^3
\right]\la_1+o(1).
\end{eqnarray*}

The   remaining task is to plug in the expansions for the $\lambda_t$'s 
given in Table~\ref{tab:lambdas}, 
simplify, and apply~\eqn{forkappa}. Since $p=o(n^{-7/11})$ we  approximate   $\lambda_1$ by  $ \frac{1}{6}n^3p^3-\frac{1}{2}n^2p^3$, whilst   for $\lambda_t$, $t\geq 2$ only the first order term is important: $\lambda_2\sim\frac{1}{4}n^4p^5 $ etc. 
 This  determines $P_{ 0,\kappa}$ and hence the coefficients in the statement of Theorem~\ref{t:main}, resulting  in the statement of Theorem~\ref{t:Gnp}  for $p<n^{-7/11-\eps}$. To relax this to  $p=o(n^{-7/11})$, we only need to note that, from this conclusion, all other terms in the series $F$ in Theorem~\ref{t:main} must  have $i_\ell/j_\ell\le 7/11$. Such terms tend to 0 for $p=o(n^{-7/11})$, and the theorem follows. 
 \qed

\bigskip

\noindent{\bf Proof of Theorem~\ref{t:Gnm}}
\smallskip

Here we extend the previous proof to obtain the probablity that $\Gnm$ contains no copies of $K_3$. The starting point of our analysis is that for $p=d/\xi$ as given by~\eqn{mudef2}, from~\eqn{rewrite} and~\eqn{X0kcond2} we have
\bel{start}
\BP(X=0\mid Y=m)\sim\frac{\BP(X=0)}{\sqrt{2\pi m}\,\BP(Y=m)}.
\ee
We need to find the asymptotics of
$\BP(X=0)$ and $\BP(Y=m)$ in the way described in the last paragraph of Section~\ref{s:maingnm}.
This requires first finding the asymptotic expansion of $p=d/\xi$, where
$\xi = \xi_{0,\epsilon}(n,p,\tilde\bfg)$.

Table~\ref{tab:cs2}
  is essentially an extension of Table~\ref{tab:cs}, showing significant contributions to \ $F_t$  from (\ref{Ftdef}) as needed to
calculate $\hatma_0(n,p,\tilde{\bfg})$, under the same assumption that $p=O(n^{-7/11-\eps})$. Note that $F_1$ and $F_2$  need  to be recomputed in this new clustering as the expression  
for $\hatma_0$ contains
$\hatma_1$ and $\hatma_2$.

Since 
$\tilde g_0=m/\la_0 =  mp^{-1} {n \choose 2}^{-1}= d/p=\xi= \xi_{0,\epsilon}(n,p,\tilde\bfg)=\hatma_0(m\delta_0) +o(\lambda_0^{-1})$ by~\eqn{hatma0approx}, 
it is straightforward to see that   the factors
$\tilde{g}_1/\hatma_0$ can at this point be replaced by $1$. Strictly this needs to be justified  in the context of the recursive computation of $\xi$ in Corollary~\ref{expanded}, and this can be seen in a straightforward way by going back to the original equations in Proposition~\ref{basic} with the altered equations and observe that the same argument as in Section~\ref{s:recursions} applies to these altered equations, resulting in the modified definition of $F_t$ in~\eqn{Ftdef}; alternatively, one could include the factors explicitly and watch them turn naturally into 1.
Note that terms like $c(1,0,  \delta_1)$ cannot affect this computation since they contain a factor $g_t$ for $t>0$, and to evaluate $\xi$ we must set such $g_t$ equal to 0.

 The denominator of (\ref{Ftdef}) is $c(0,0,\bz)=1-p$ in the case of $t=0$. 
As with the $\Gnp$ calculation, we can ignore certain  terms in the product of $c(u,t,h)$ with its cofactor. In the case of $t=0$, since the final expression we are computing is $\exp(-\la_1+o(\la_1))$ and $\la_1\sim n^3p^3/6$,    we can   ignore any terms that are   $O(\error/\la_1)$, i.e.\ $O(\error/n^3p^3)$. Note that $\la_0=n(n-1)p/2$.
Since $\hatma_1$ only  arises in terms with a cofactor that is $O(np^2)$, we ignore terms in its expression that are $O(\error/n^4p^5)$ such as $p^2$. For similar reasons, terms in $\hatma_2$ of order  $O(\error/n^5p^7)$ are ignored.

\def\fracspace{\raisebox{-2mm}{\rule{0cm}{6mm}}}
 \def\exspace{\fracspace}
\begin{table}[h]
 \small
\begin{center}
\begin{tabular}{|c|c|c|c|c|}
\hline
\exspace  $u$& t&h&$c(u,t,h)$& cofactor\\
\hline
\fracspace 0&0&$\delta_0$& $ p $& $p \cdot \hatma_0$\\
\hline
\fracspace 1&0&$2\delta_0$& $ 3(1-p)$& $\frac13  np^2 \cdot  \hatma_1 $  \\
\hline
\exspace 2&0&$4\delta_0$& $ 1 $& $\frac{1}{2}n^2p^4 \cdot  \hatma_2 $\\
\hline
\exspace 4&0&$6\delta_0$& $ 1 $& $\frac{1}{6}n^3p^6$  \\
\hline
\fracspace  2&1&$2\delta_0$& $ 2 $& $\frac32 np^2 \cdot   \hatma_2 $\\
\hline
\exspace 3&1&$4\delta_0$& $ 1 $& $3n^2p^4$ \\
\hline
\exspace 4&1&$4\delta_0$&$3$    &$\frac12 n^2p^4 $ \\ 
\hline
\exspace 3&2&$2\delta_0$      & $2$ & $2np^2 $\\
\hline
\exspace 4&2&$4\delta_0$& $3$ & $\frac13 np^2  $\\
\hline
\end{tabular}

\end{center}
\caption{Significant contributions to~\eqn{Ftdef}}
\lab{tab:cs2}
\end{table}

Plugging  the values in Table~\ref{tab:cs2} 
 into (\ref{hatmadef}) or~\eqn{Ftdef}
gives the following truncated expressions 
for $\hatma_0$, $\hatma_1$ and $\hatma_2$ 
$$
\hatma_0=\frac{1}{1-p}\left(1-p\hatma_0-(1-p)np^2\hatma_1
-\frac{1}{2}n^2p^4\hatma_2-\frac{1}{6}n^3p^6
\right),
$$
$$
\hatma_1=1-3np^2\hatma_2 -\frac{9}{2}n^2p^4,
$$
and
$$
\hatma_2=1-5np^2.
$$
Solving for $\hatma_0$ gives
\bean
\hatma_0 = 1-np^2-\frac{49}{6}n^3p^6+\frac{5}{2}n^2p^4+np^3
+\frac{21}{2}n^3p^7-3n^2p^5.
\eean
Using this  expression for $\xi$ we find $1/\xi\approx  1+x-\frac32 x^2+\frac{25}{6}x^3-x p$ where $x=np^2$ and the terms of order $x^4$, $x^2p$ and $p^2$ are omitted. Substituting $p=d/\xi$ into itself three times
gives
\eq\lab{pnmtri}
p = d+nd^3-nd^4+\frac{1}{2}n^2d^5+\frac{1}{6}n^3d^7 +o(1/n^3p^3).
\en
  
For $N={{n}\choose{2}}$ and noting $p=d(1+\eps)=m(1+\eps)/N$ where $\epsilon=O(d^2n)$, we have
\begin{eqnarray}
 \BP(Y=m)&=& \BP(Y=dN) = {N\choose m} 
p^{dN}{(1-p)}^{ N(1-d)} \\
&\sim &
\frac{1}{\sqrt{2\pi dN}}\left((1+\eps)^d\left(\frac{1-d(1+\eps)}{1-d}\right)^{1-d}\right)^N
\nonumber\\
&\sim&
\frac{1}{\sqrt{2\pi m}}
\exp\left( 
-\frac14n^4 d^5-\frac{1}{12}n^5 d^7-\frac{1}{48}n^6 d^9+\frac14n^4d^6
\right)  \lab{prk}
\end{eqnarray}
using~\eqn{pnmtri} to determine $\eps$.
 Theorem \ref{t:Gnp} gives $\BP(X=0)$, from which we can again eliminate $p$ using~\eqn{pnmtri}.  Plugging these into~\eqn{start}  gives the probability that $X=0$ in $G\in{\cal G}(m,n)$:
\[
\BP(X=0|Y=m)\sim
\exp\left(
-\frac{1}{6}n^3d^3  -\frac{1}{8}n^4d^6\right).
\]
For the same reasons as in the $\Gnp$ case, the validity extends to all $d=o(n^{-7/11})$.  \qed

\end{appendices}


\begin{thebibliography}{99}


\bibitem{B} B.~Bollob{\'a}s, Random graphs, In {\em Combinatorics}, Proceedings (Swansea, 1981), pp. 80--102, London Math.\ Soc.\ Lecture Note Ser.\ 52, Cambridge Univ.\ Press, Cambridge, 1981.


\bibitem{EKR} P. Erd{\H o}s, D.J. Kleitman and B.L. Rothschild, Asymptotic enumeration of $K_n$-free graphs, Colloquio Internazionale sulle Teorie Combinatorie (Rome, 1973), Tomo II, pp. 19Ð27. Atti dei Convegni Lincei, No. 17, Accad. Naz. Lincei, Rome, 1976.

\bibitem{F} A.~Frieze, On small subgraphs of random graphs. In
{\em Random Graphs, Volume 2}, Wiley, New Tork, (1992), 67--90.


\bibitem{JLRinequal} S.~Janson, T. \L uczak and A. Ruci\'nski,
An exponential bound for the probability of nonexistence of a specified subgraph in a random graph, in: M. Karo{\'n}ski, J. Jaworski, A. Ruci{\'n}ski, eds., Random Graphs '87 (Wiley), pp. 73--87, 1990.

\bibitem{JLR} S.~Janson, T. \L uczak and A. Ruci\'nski,
{\em Random graphs}. Wiley, New York, 2000.



\bibitem{L} T. {\L}uczak, On triangle-free random graphs, {\em Random Structures \& Algorithms} {\bf 16} (2000), 260--276.

\bibitem{OPT}  D. Osthus, H.J.~Pr\"omel and A. Taraz,  For which densities are random triangle-free graphs almost surely bipartite? {\em Paul Erd{\H o}s and his mathematics (Budapest, 1999)}. {\em Combinatorica} {\bf 23},  105--150. 

\bibitem{PS} H.J.~Pr\"omel and A.~Steger~\cite{PS}
Counting $H$-free graphs.  
{\em Discrete Math.}, {\bf 154}  (1996), 311--315.

\bibitem{PS96} H.J.~Pr\"omel and A.~Steger,
On the asymptotic structure of sparse triangle free graphs,
{\em J. Graph Theory}, {\bf 21} (1996), 137--151.

 \bibitem{R}   A. Ruci{\'n}ski, When are small subgraphs of a random graph normally distributed? {\em Probab.\ Theory Related Fields} {\bf 78} (1988),  1--10.

\bibitem{W} N.C. Wormald, The perturbation method and triangle-free random
graphs, {\em Random Structures \& Algorithms}, {\bf 9} (1996), 253--270.

\end{thebibliography}
\end{document}